\documentclass[10pt,a4paper,reqno]{amsart} 
\usepackage{latexsym}
\usepackage{gensymb}
\usepackage{verbatim}
\usepackage{epsfig}
\usepackage{rotating}
\usepackage{amssymb}
\usepackage[T1]{fontenc}
\usepackage{kpfonts}
\usepackage{afterpage}
\usepackage{color}
\usepackage{dsfont}
\usepackage{amssymb,amsfonts,amsmath,stmaryrd,bbm}
\usepackage[utf8]{inputenc}
\usepackage{pifont, tikz, subfigure}
\usetikzlibrary{plotmarks,shapes,arrows,positioning}
\usepackage{url}
\usepackage{pifont, tikz, subfigure}
\usetikzlibrary{plotmarks,shapes,arrows,positioning}

\usepackage[plainpages=false,pdfpagelabels,colorlinks=true,citecolor=blue,hypertexnames=false]{hyperref}

\newcommand{\tS}{\tilde S}

\catcode`\@=11
\def\section{\@startsection{section}{1}%
 \z@{.7\linespacing\@plus\linespacing}{.5\linespacing}%
 {\normalfont\Large\bfseries\scshape\centering}}

\def\subsection{\@startsection{subsection}{2}%
  \z@{.5\linespacing\@plus\linespacing}{.5\linespacing}%
  {\normalfont\large\bfseries\scshape}}

\def\subsubsection{\@startsection{subsubsection}{3}%
 \z@{.5\linespacing\@plus\linespacing}{-.5em}
 {\normalfont\large\bfseries}}
\catcode`\@=12

\newtheorem{theorem}{Theorem}

\newtheorem{cor}[theorem]{Corollary}
\newtheorem{prop}[theorem]{Proposition}

\theoremstyle{definition}

\newcommand{\beq}{\begin{equation}}
\newcommand{\eeq}{\end{equation}}

\renewcommand{\epsilon}{\varepsilon}
\renewcommand{\iota}{\phi}

\makeatletter
\def\testb#1{\testb@i#1,,\@nil}%
\def\testb@i#1,#2,#3\@nil{%
  \draw[->, thick] (O) --++(#1);
  \ifx\relax#2\relax\else\testb@i#2,#3\@nil\fi}
\makeatother   
\newcommand{\makediag}[1]{
    \coordinate (O) at (0,0); \coordinate (N) at (0,0.8);
    \coordinate (NE) at (0.8,0.8); \coordinate (E) at (0.8,0);
    \coordinate (SE) at (0.8,-0.8); \coordinate (S) at (0,-0.8);
    \coordinate (SW) at (-0.8,-0.8);\coordinate (W) at (-0.8,0);
    \coordinate (NW) at (-0.8,0.8); \coordinate (B1) at (1.2,1.2);
    \coordinate (B2) at (-1.2,-1.2);
    \testb{#1}
} 
\newcommand{\diagr}[1]{
  \begin{tikzpicture}[scale=0.8]\makediag{#1}\end{tikzpicture}
}

\def\sequenceThreeD#1#2#3#4{%
  #2\quad\hbox to.6\hsize{$#3,\dots$\hfill}\quad\textrm{(#4)}
}

\long\def\greybox#1{%
    \newbox\contentbox%
    \newbox\bkgdbox%
    \setbox\contentbox\hbox to \hsize{%
        \vtop{
            \kern\columnsep
            \hbox to \hsize{%
                \kern\columnsep%
                \advance\hsize by -2\columnsep%
                \setlength{\textwidth}{\hsize}%
                \vbox{
                    \parindent=0bp
                    #1
                }%
                \kern\columnsep%
            }%
            \kern\columnsep%
        }%
    }%
    \setbox\bkgdbox\vbox{
        \pdfliteral{0.9 0.9 0.9 rg}
        \hrule width  \wd\contentbox %
               height \ht\contentbox %
               depth  \dp\contentbox
        \pdfliteral{0 0 0 rg}
    }%
    \wd\bkgdbox=0bp%
    \vbox{\hbox to \hsize{\box\bkgdbox\box\contentbox}}%
}

\newcommand{\zs}{\mathbb{Z}}

\newcommand{\qs}{\mathbb{Q}}

\newcommand{\cs}{\mathbb{C}}
\newcommand{\fps}{formal power series}

\newcommand{\bx}{\bar x}
\newcommand{\by}{\bar y}

\newcommand{\Sp}{S_1}
\newcommand{\Ch}{C_{-0}}
\newcommand{\Cv}{C_{0-}}
\newcommand{\Ah}{A_{-0}}
\newcommand{\Av}{A_{0-}}

\newcommand{\cC}{\mathcal C}
\newcommand{\cS}{\mathcal S}

\DeclareMathOperator{\Pol}{Pol}
\DeclareMathOperator{\Rat}{Rat}

\newcommand{\gf}{generating function}
\newcommand{\gfs}{generating functions}

\def\emm#1,{{\em #1}}

\newcommand{\C}{C}
\newcommand{\Cn}{C_-}
\newcommand{\Qp}{Q_+}

\begin{document}
%
\title[Square lattice walks avoiding a quadrant]
{Square lattice walks avoiding a quadrant}

\author[M. Bousquet-M\'elou]{Mireille Bousquet-M\'elou}

\thanks{}
 
\address{CNRS, LaBRI, Universit\'e de Bordeaux, 351 cours de la
  Lib\'eration,  F-33405 Talence Cedex, France} 
\email{bousquet@labri.fr}

\begin{abstract} In the past decade, a lot of attention has been
  devoted to the enumeration of walks with prescribed steps confined
  to a convex cone. In two dimensions, this means  counting walks in
  the first quadrant of the plane (possibly after a
  linear transformation).

But what about walks in non-convex cones?  We
investigate the two most natural cases: first, square lattice walks avoiding
the negative quadrant $\mathcal Q_1=\{(i,j): i<0 \hbox{ and } j<0\}$, and then,
square lattice walks avoiding the West quadrant $\mathcal Q_2=\{(i,j): i <j \hbox{
  and } i<-j\}$.  In both cases,  the  generating function
that counts walks starting from the origin 
is found to differ from a  simple D-finite series by an
algebraic one. 
We also obtain  closed form
expressions for the number of $n$-step walks ending at certain prescribed endpoints,  as a sum of
three hypergeometric terms. 

One of these terms already appears in the enumeration of square
lattice walks confined to the cone $\{(i,j): i+j \ge 0 \hbox{ and }
j\ge 0\}$, known as Gessel's walks. In fact, the enumeration of
Gessel's walks follows, by the reflection principle, from the enumeration of walks starting from
$(-1,0)$ and avoiding $\mathcal Q_1$. Their generating function turns
out to be purely algebraic (as the \gf\ of Gessel's walks).

Another approach to Gessel's walks consists in counting walks that start
from $(-1,1)$ and avoid the West quadrant $\mathcal Q_2$. The
associated \gf\ is D-finite but transcendental.
\end{abstract}

\maketitle

\section{Introduction}
\begin{figure}[hbt]
\begin{center}
\includegraphics{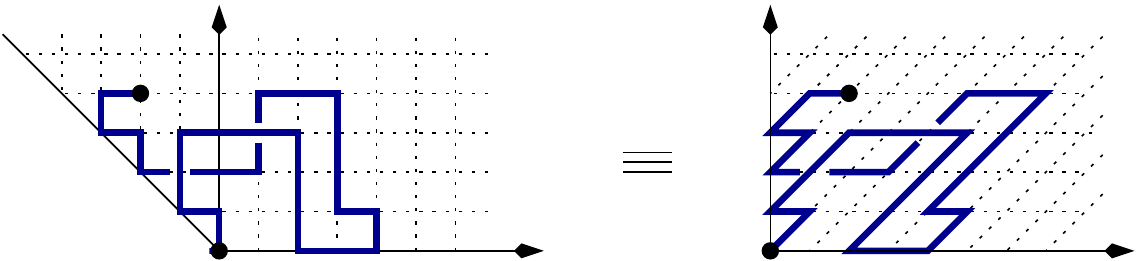}
\end{center}
\caption{Square lattice walks staying in a $135\degree$ wedge are equivalent
  to quadrant walks with steps $\rightarrow, \nearrow, \leftarrow, \swarrow$.}
\label{fig:gessel}
\end{figure}

In recent years, the enumeration of lattice walks confined to convex
cones has attracted a lot of attention.  In two dimensions, this means
counting walks in the intersection of two half-spaces, which we can
always assume
(Figure~\ref{fig:gessel}) to form the first quadrant $
\mathcal Q=\{(i,j): i\ge 0 \hbox{ and } j \ge 0\}$. 
 The problem is then completely specified by
prescribing a starting point and a set of allowed
steps. The  two most natural examples are walks on the  square lattice
(with
steps $\rightarrow, \uparrow, \leftarrow, \downarrow$), and walks on the
\emm diagonal, square lattice (with steps $\nearrow$, $\nwarrow$,
$\swarrow$, $\searrow$). Both cases can be solved via the classical
reflection principle~\cite{gessel-zeilberger,guy-bijections}. The enumeration usually records the length $n$ of
the walk (with a variable $t$), and the coordinates $(i,j)$ of its endpoint
(with variables $x$ and $y$). For instance, the \gf\
of square lattice walks starting from $(0,0)$ and confined to
$\mathcal Q$ is~\cite{guy-bijections,BoMi10}:
 \beq\label{Q-expr}
Q(x,y)=\sum_{i,j,n\ge 0}
\frac{(i+1)(j+1)}{(n+1)(n+2)}\binom{n+2}{\frac{n-i-j}2}\binom{n+2}{ \frac{n+i-j+2}2}
x^i y^j t^{n},
\eeq
where the sum is restricted to integers $i,j, n$ such that $n$ and
$i+j$ have the same parity.  (To
lighten notation, we ignore the dependence in $t$ of this
series.) This  series  is \emm D-finite,~\cite{lipshitz-df}: this means that it satisfies a linear
differential equation in each of its variables $t$, $x$ and $y$, with
 coefficients in the field $\qs(t,x,y)$ of rational functions in $t$,
 $x$ and $y$.

In the past decade, a systematic study of quadrant walks with \emm
small steps, (that is, steps in $\{-1,0,1\}^2$) has been carried out,
and a complete classification is now available. For walks starting at
$(0,0)$, the \gf\ is D-finite if and only if a certain group of
birational transformations is finite. The proof combines an attractive combination of
approaches:  algebraic~\cite{bousquet-versailles,BoMi10,gessel-proba,gessel-zeilberger,mishna-jcta,niederhausen-ballot},
computer-algebraic~\cite{BoKa10,KaKoZe08,kauers07v},
analytic~\cite{BoKuRa13,KuRa12,raschel-unified},
asymptotic~\cite{BoRaSa12,denisov-wachtel,MeMi13,MiRe09}.

The most intriguing D-finite case is probably Gessel's model,
illustrated in Figure~\ref{fig:gessel}. Around 2000, Ira Gessel conjectured that
the number of $2n$-step walks of this type starting and ending at $(0,0)$ was
\beq\label{conj:gessel}
g_{0,0}(2n)=16^n\, \frac{ (1/2)_n(5/6)_n}{(2)_n(5/3)_n},
\eeq
where $(a)_n=a(a+1) \cdots (a+n-1)$ is the ascending factorial.
A computer-aided proof of this 
conjecture was finally found in 2009 by Kauers, Koutschan and
Zeilberger~\cite{KaKoZe08}. A year later, Bostan and
Kauers~\cite{BoKa09} proved, using again intensive computer algebra, that the 
three-variate \gf\ of Gessel's walks starting at $(0,0)$ and ending
anywhere in the quadrant
is not only D-finite, but even algebraic: this means that it satisfies a
polynomial equation over $\qs(t,x,y)$. Three other proofs have now been
given~\cite{BeBoRa16,BoKuRa13,mbm-gessel}, but none of them explains
combinatorially the
simplicity of the numbers, nor the algebraicity of the series.

\medskip

\begin{figure}[b!]
\begin{center}
\includegraphics{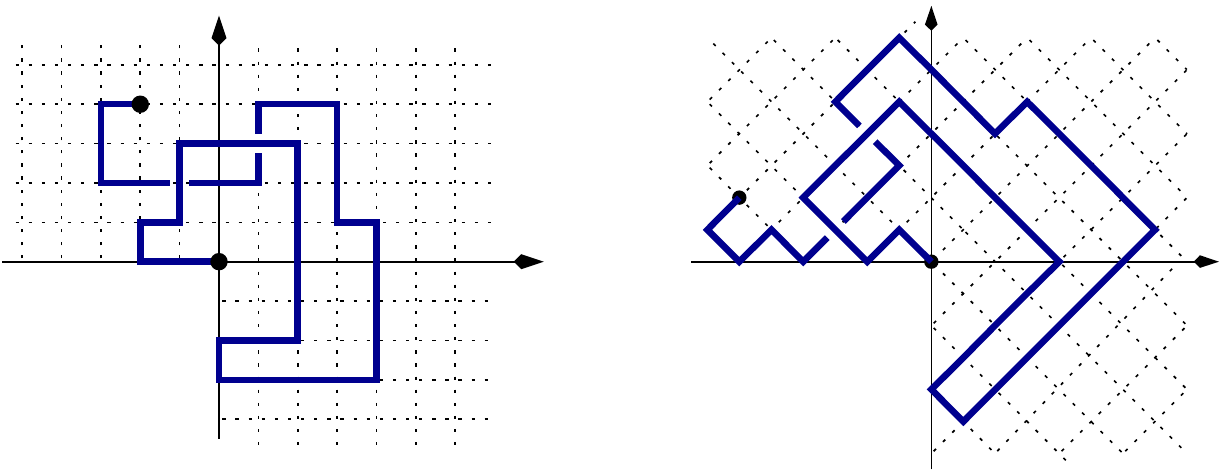}
\end{center}
\caption{Walks  confined to the non-convex cone $\mathcal \C$: on the
  square lattice (left), and on the diagonal square lattice (right).}
\label{fig:square-diag}
\end{figure}

The primary objective of this paper is to initiate a parallel study
for walks confined to non-convex cones. In two dimensions, this means
that walks live in the \emm union, of two half-spaces, which we can
assume to form   the three-quadrant cone 
$$
\mathcal \C:=\{(i,j): i\ge 0 \hbox{ or } j \ge 0\}.
$$
In other words, these walks avoid the negative quadrant. We solve here
the two most natural cases (and possibly the simplest): the square lattice, and the diagonal
square lattice (Figure~\ref{fig:square-diag}). By a simple rotation,
 the latter model is equivalent to square lattice walks
avoiding the West quadrant $
\{(i,j) : i < -|j|\},
$
as described in the abstract. The first problem was raised in 2001
by David W. Wilson in entry A060898 of the OEIS~\cite{oeis}. 

These two
problems are far from being as simple as their quadrant
counterparts. Their solutions exhibit, as in Gessel's
model, a combination of algebraicity phenomena and  
hypergeometric numbers. For instance, the \gf\ $C(x,y)$ of square
lattice walks starting at $(0,0)$ and confined to $\cC$ differs from
the D-finite series
$$
\frac 1 3 \left(Q(x,y) -\bx^2 Q(\bx,y) - \by^2
Q(x,\by) \right)
$$
by an \emm algebraic, one (we have written $\bx$ for $1/x$ and $\by$ for
$1/y$, and the series $Q$ is given by~\eqref{Q-expr})). This holds as well on the diagonal square lattice, if $Q(x,y)$ now
counts quadrant walks with diagonal steps. In terms of numbers, we
find for instance that the number of walks of length $2n$ starting and ending at
$(0,0)$ on the diagonal square lattice is 
 $$
c_{0,0}(2n)=\frac{16^n}9 \left(  3 \, \frac{(1/2)_n^2}{(2)_n^2}
+  8 \,  \frac{(1/2)_n(7/6)_n}{(2)_n(4/3)_n}
-  2 \,  \frac{(1/2)_n(5/6)_n}{(2)_n(5/3)_n}\right).
$$
Two ingredients of this formula are familiar: 
 $16^n\frac{(1/2)_n^2}{(2)_n^2}$ counts walks confined to the
  first quadrant, while $16^n\, \frac{ (1/2)_n(5/6)_n}{(2)_n(5/3)_n
    }$ counts
 Gessel's walks.
Asymptotically, these two terms are dominated by the central one,
and
$$
c_{0,0}(2n)\sim \frac{2^5}{3^2} \frac {\Gamma(2/3)}{\pi} \frac {4^{2n}}{(2n)^{5/3}}.
$$
We obtain a similar, slightly more
complicated formula for  the square lattice (see~\eqref{c00-square}).

That Gessel's numbers are involved in this problem should not be
too surprising. Indeed, halving three quadrants gives a $135\degree$ cone as in
Figure~\ref{fig:gessel},  and the solution of Gessel's problem can be recovered from the reflection principle  if we count square lattice walks  starting from
$(-1,0)$ and confined to $\cC$ (Figure~\ref{fig:refl}). An alternative approach is to count walks on the
\emm diagonal, square lattice starting from $(-2,0)$ and confined to
$\cC$ (Figure~\ref{fig:refl-diag}). This connection between
three-quadrant problems and Gessel's walks was in fact another 
motivation of our study, and we solve these two problems with shifted
starting point. We do not claim to have explained
combinatorially Gessel's ex-conjecture:  as all proofs of this
conjecture, our approach to the three-quadrant problem consists in
solving a functional equation satisfied by the \gf\ $\C(x,y)$. The
tools involved in the solution consist of elementary power series
manipulations, coefficient extractions, polynomial elimination. We have
at the moment no 
combinatorial understanding of our results.

\medskip
 We hope that this work will be the starting point of a systematic 
study of walks avoiding a quadrant, analogous to what has been done so far
for walks confined to a quadrant.  The difficulty of the
``simple'' square lattice case 
suggests that this study may turn out to be even
more challenging.

The paper is organized as follows. 
In Section~\ref{sec:square} we count square lattice walks starting
from $(0,0)$ and confined to $\cC$. Analogous
results are proved in Section~\ref{sec:diag} for
walks on the diagonal square lattice. 
In  Section~\ref{sec:square-asym} we go back to the square lattice,
but change the starting point to $(-1,0)$.  The $x/y$ symmetry is
lost, but we still obtain a complete solution, in fact simpler 
than in the original case:
\begin{quote}
  \emph{the \gf\ of square lattice walks that start from $(-1,0)$ and avoid
  the negative quadrant is algebraic.}
\end{quote}
In  Section~\ref{sec:diag-asym} we count similarly walks starting from
$(-2,0)$ on the diagonal square lattice. The \gf\ is not algebraic,
but differs from a simple D-finite series by an algebraic one.  In Section~\ref{sec:gessel} we derive from
these results a new solution of Gessel's problem. 
Some perspectives and open questions are
discussed in the final section. The paper is  accompanied by two
Maple sessions (one for the square lattice, one for the diagonal square
lattice)
  available on the author's \href{http://www.labri.fr/perso/bousquet/publis.html}{webpage}.

\medskip
\noindent {\bf Notation.} 
For a ring $R$, we denote by $R[t]$
 (resp.~$R[[t]]$) 
the ring of
polynomials (resp. formal power series) 
in $t$ with coefficients in
$R$. If $R$ is a field, then $R(t)$ stands for the field of rational functions
in $t$. 
This notation is generalized to several variables.
For instance, the \gf\ $Q(x,y)$ that counts quadrant walks is a series
of $\qs[x,y][[t]]$, while the \gf\ 
$C(x,y)$ of walks confined to~$\cC$ belongs to $\qs[x,\bx,y,\by ][[t]]$, where $\bx=1/x$ and
$\by=1/y$.

If $G(x)$  is a power series in~$t$ with coefficients in $\qs[x,
  \bx]$, written as
$$
G(x)=\sum_{i\in \zs, \ n\ge 0}g_i(n) t^n x^i,
$$
 we denote by $G_i$ the
coefficient of $x^i$ in $G$:
\beq\label{Gi-def}
G_i= [x^i]G(x):=  \sum_{n\ge 0}  g_i(n) t^n.
\eeq
We denote by $[x^{>}]G(x)$ the \emm positive part of $G$ in $x$,:
$$
 [x^{>}]G(x):= \sum_{i>0, \ n\ge 0}  g_i(n) t^n x^i.
$$
We similarly  define the non-positive, negative, and non-negative parts
of $G(x)$. Finally, for a series
$G(x,y)$ in $t$ with  coefficients in $\qs[x,\bx,y,\by]$, we
denote by $G_{i,j}$ the coefficient of $x^i y^j$, which is a series in
$t$. 

We refer to~\cite{lipshitz-df} for properties of D-finite series.

\section{The square lattice}
\label{sec:square}
The aim of this section is to determine the \gf\ of square lattice
walks starting from $(0,0)$ and confined to the three-quadrant cone $\cC$. It reads
$$
\C(x,y)=  \sum _{(i,j) \in \mathcal \C}\sum_{n \ge 0 }c_{i,j}(n) x^i
y^j t^n
= 1+t (x+\bx+y+\by) + O(t^2),
$$
where $c_{i,j}(n)$ counts $n$-step walks going from $(0,0)$ to
$(i,j)$. For walks confined to the first quadrant $\mathcal Q$, we define
the series $Q(x,y)$ and its coefficients $q_{i,j}(n)$
similarly. As recalled in the introduction,  these
coefficients have a simple hypergeometric form (see~\eqref{Q-expr}).

 Our first result
(Theorem~\ref{thm:square}) 
states that the \gf \ $C(x,y)$  differs from
the  simple D-finite series 
$$
\frac 1 3 \left(Q(x,y) -\bx^2 Q(\bx,y) - \by^2
Q(x,\by) \right)
$$
by an  algebraic series, which we describe explicitly. From there, we
can express the \gf\ $C_{i,j}$ of
walks ending at a prescribed point $(i,j)$ (Corollary~\ref{cor:coeffs}), and in some cases, obtain closed
form expressions for its coefficients $c_{i,j}(n)$. 

\begin{theorem}\label{thm:square}
The \gf\ of square lattice walks starting at $(0,0)$, confined to $\mathcal \C$ 
and ending in the first quadrant
(resp. at a negative abscissa) is
\beq\label{C-split}
\frac 1 3\,Q(x,y) + P(x,y), \qquad\qquad
\left( \hbox{resp.} \quad -\frac 1 3\, \bx^2 Q(\bx,y)+  \bx M(\bx,y)
\right),
\eeq
where $M(x,y)$ and $P(x,y)$ are algebraic  of degree $72$ over $\qs(x,y,t)$.

More precisely,  $P$ can be expressed in terms of $M$ by:
\beq\label{PM}
P(x,y)=\bx\left( M(x,y)-M(0,y)\right) + \by \left(
M(y,x)-M(0,x)\right),
\eeq
$M$ satisfies the functional equation
\begin{multline}\label{func-M}
  (1-t(x+\bx+y+\by)) \left( 2M(x,y)-M(0,y)\right) =\\ 2x/3 -2t \by M(x,0) + t(x-\bx)
M(0,y)  +t\by M(y,0),
\end{multline}
and the specializations $M(x,0)$ and $M(0,y)$ have respective degrees
$24$ and $12$  over $\qs(t,x)$ and $\qs(t,y)$. 

Moreover, these algebraic  series admit
rational parametrizations. Let $T$ be the unique series in $t$ with
constant term $1$ satisfying:
\beq\label{T-def}
T= 1+256\, t^2 \frac{T^3}{(T+3)^3},
\eeq
and let $Z=\sqrt T$.
Let  $U$ be the only power series in $t$ with constant term $1$ 
satisfying 
\beq\label{param-x}
16T^2(U^2-T)=
x(U+UT-2T)(U^2-9T+8TU+T^2-TU^2).
\eeq
Then the series $tM(xt,0)$ and $tM(0,xt)$ (both even series in $t$)
admit rational expressions in terms of $Z$ and $U$, given in  Appendix~\ref{sec:param-sq}.
\end{theorem}

\noindent{\bf Remarks}\\
1. Equation~\eqref{C-split} gives the \gfs\ of walks ending in two quadrants
of $\cC$. By symmetry, the \gf\ of walks ending in the third quadrant,
that is, at a negative ordinate, is $-\by^2Q(\by,x)/3+\by M(\by,x)$.
\\
2. The parametrization by $T$, $Z$ and $U$ already
appears in van Hoeij's parametrization for Gessel's walks in the
quadrant~\cite{BoKa10,mbm-gessel}. As explained later in
Section~\ref{sec:gessel}, this is no coincidence. A fourth series is
involved in the parametrization of Gessel's problem, and we will find
it when counting walks confined to $\cC$ on the diagonal square lattice (Theorem~\ref{thm:diag}).\\
3. A more compact statement of Theorem~\ref{thm:square} reads as follows:
$$
C(x,y)= A(x,y) + \frac 1 3 \left(Q(x,y) -\bx^2 Q(\bx,y) - \by^2
Q(x,\by) \right),
$$
where $A(x,y)$ satisfies
$$
(1-t(x+\bx+y+\by))A(x,y)= (2+\bx^2+\by^2)/3-t\by A_-(\bx) -t \bx A_-(\by) ,
$$
and $A_-(x)$ is a series in $t$ with coefficients in $\qs[x]$,
algebraic of degree 24. It equals the series $xM(x,0)$ given in
Appendix~\ref{sec:param-sq}. 

Of course  the algebraicity of $A(x,y)=P(x,y)+\bx M(\bx,y)+\by
M(\by,x)$ follows from this statement, but it  hides the fact that
the series $P(x,y)$ and $M(x,y)$ are algebraic themselves (there is no
reason why extracting say, negative powers of $x$ in an algebraic
series with coefficients in $\qs[x,\bx,y,\by]$ should yield an
algebraic series). \\
4. The four series $Q(x,y), Q(\bx,y), Q(x,\by)$ and $Q(\bx,\by)$ are
related by a simple identity (see~\eqref{Q-lin}), which allows us to write:
$$
C(x,y)=A(x,y)- \frac 1 3 \bx^2 \by ^2 Q(\bx,\by)+
\frac{(x-\bx)(y-\by)}{3xy (1-t(x+\bx+y+\by))}.
$$
 This implies that $C(x,y)$,
as $Q(x,y)$ itself, is D-finite but transcendental.
 \\
5. By combining the above results and \emm singularity analysis, of
algebraic (and D-finite) series~\cite{flajolet-sedgewick}, one can
derive from the above theorem that the number
of $n$-step walks confined to $\cC$ is
$$
[t^n] C(1,1)\sim \frac{2^5 \sqrt 3}{3^3\, \Gamma(2/3)}
\frac{4^n}{n^{1/3}}.
$$ 

\medskip
We now focus on walks ending at a prescribed point. For $(i,j) \in
\mathcal C$, let $\C_{i,j}$
denote the length \gf\ of walks going
from $(0,0)$ to $(i,j)$  in
$\mathcal \C$. Define similarly $Q_{i,j}$, for $i, j \ge 0$. According
to~\eqref{Q-expr}, the latter series has a simple form:
$$
Q_{i,j}=\sum_{n\ge 0}
\frac{(i+1)(j+1)}{(n+1)(n+2)}\binom{n+2}{\frac{n-i-j}2}\binom{n+2}{ \frac{n+i-j+2}2}
 t^{n}.
$$
The following corollary clarifies the
nature of the series $\C_{i,j}$.

\begin{cor}[\bf Walks ending at a prescribed position]\label{cor:coeffs}
  Let $T$ be the unique series in $t$ with constant term $1$
  satisfying~\eqref{T-def}, and let  $Z=\sqrt T $.

 For $j\ge 0$, the series $C_{-1,j}$ belongs  to $t^{j+1} \qs(Z)$, and
 is thus algebraic.
More generally, for $i\ge 1$ and $j\ge 0$, the series $C_{-i,j}$ is
D-finite, of the form
 $$
- \frac 1 3 Q_{i-2,j} + t^{i+j} \Rat(Z)
$$
for some rational function $\Rat$. It is transcendental
as soon as $i\ge 2$.

Finally, for $i\ge 0$ and $j\ge 0$,  the series $C_{i,j}$  is of the form
 $$
 \frac 1 3 Q_{i,j} + t^{i+j} \Rat(Z).
$$
It is D-finite and transcendental.
\end{cor}
\noindent The series $C_{i,j}$ can be effectively computed. For instance,
\beq\label{coincidence1}
tC_{-1,0}=
 {\frac {(Z^2-1) \left( 11+6 {Z}^{2}-{Z}^{4} \right) }{ \left( {Z}^{2}+3 \right) ^{3}}}
,
\eeq
$$
C_{-1,1}=
1024 {\frac {{Z}^{3} \left( {Z}^{2}+1 \right) ^{2} \left( Z-1 \right)  
\left(1+2Z- {Z}^{2} \right) }{ \left( {Z}^{2}+3 \right) ^{6} \left( Z+1 \right) }},
$$
$$
C_{-2,0}= -\frac 1  3 Q_{0,0}+{\frac {256\, {Z}^{3} 
\left(4+4 Z-4
  {Z}^{2} +23 {Z}^{3}-9 {Z}^{4}+18 {Z}^{5}-6 {Z}^{6} +3 {Z}^{7}-{Z}^{8}\right) }
{3  \left( {Z}^{2}+3 \right) ^{6} \left( Z+1 \right) }},
$$
$$
C_ {0,0}= \frac 1 3 Q_{0,0} +{\frac {512\, {Z}^{3} 
\left(4+4 Z-4
  {Z}^{2} +23 {Z}^{3}-9 {Z}^{4}+18 {Z}^{5}-6 {Z}^{6} +3 {Z}^{7}-{Z}^{8}\right) 
}{3  \left( {Z}^{2}+3 \right) ^{6} \left( Z+1 \right) }}.
$$
The similarity between the last two expressions comes from~\eqref{PM},
which  tells us that $P_{0,0}=2M_{1,0}$, while $C_{0,0}= Q_{0,0}/3+ P_{0,0}$ and
$C_{-2,0}=-Q_{0,0}/3+ M_{1,0}$ by~\eqref{C-split}.

\medskip
Starting from the expression of $C_{i,j}$ (more precisely, of its
algebraic part $P_{i,j}$ or $M_{-i-1,j}$, depending on the sign of
$i$), one can decide if the coefficients $c_{i,j}(n)$ have an
expression as a finite sum of hypergeometric terms: one first
computes a linear recurrence relation with polynomial coefficients
satisfied by the coefficients (for instance using the Maple
commands {\tt algeqtodiffeq} and {\tt diffeqtorec}) and then applies the
Hyper algorithm from~\cite{AB}, which determines all hypergeometric
solutions of such a recurrence relation. Using    its Maple
incarnation {\tt hypergeomsols}, we  thus obtain:
\begin{align}
  c_{0,0}(2n)&= \frac{4\cdot 16^n}{3^5}\left( 3^4 \frac {(1/2)_n
    (1/2)_{n+1}}{(2)_n(2)_{n+1}} + 4 (24n^2+60n +29)
  \frac{(1/2)_n(7/6)_n}{(2)_{n+1}(4/3)_{n+1}}\right. \nonumber
\\
&\hskip 50mm \left.
-2 (12n^2+30n+5)  \frac{(1/2)_n(5/6)_n}{(2)_{n+1}(5/3)_{n+1}}\right)
\label{c00-square}
\\
&\sim \frac{2^9}{3^4} \frac{\Gamma (2/3)}{\pi} \frac {4^{2n}}{(2n)^{5/3}}.\nonumber
\end{align}
Given the link between $C_{0,0}$ and $C_{-2,0}$, this gives
 a closed form expression of the same type for
$c_{-2,0}(2n)$. We have found similar expressions for the endpoints
$(1,1)$ and $(-3,1)$, but not for $(-1,0)$, $(-3,0)$, $(-4,0)$, $(0,1)$,
$(-1,1)$, nor $(-2,1)$.

\medskip
As in the systematic study of quadrant models~\cite{BoMi10}, the starting point of
our approach is a functional equation that translates the  step by
step construction of walks
confined to $\mathcal C$. It reads:
$$
\C(x,y)=1+t(x+\bx+y+\by) \C(x,y) -t\by \Cn(\bx)-t\bx \Cn(\by),
$$
where
\beq\label{Tm-def}
\Cn(\bx)= \sum_{i<0, n\ge 0} c_{i,0}(n) x^i t^n \ \ \in \bx \qs[\bx][[t]]
\eeq
counts walks ending on the negative $x$-axis. The terms
$t\by \Cn(\bx)$ and $t\bx \Cn(\by)$ correspond to forbidden moves
yielding in the negative quadrant.
Equivalently,
\beq\label{eq-T}
K(x,y)\C(x,y)=1-t\by \Cn(\bx)-t\bx \Cn(\by)
\eeq
where
$$
K(x,y)=1-t(x+\bx+y+\by)
$$
is the \emm kernel, of the equation.

Before going further, we review a solution of the associated
quadrant model, which  adapts to many other quadrant
models~\cite{BoMi10}.  In the square lattice case which we consider here, it is
essentially a power series version of 
the classical reflection principle~\cite{gessel-zeilberger}.

\subsection{Warming up: walks confined to the positive quadrant} 
\label{sec:quadrant}
We start from the functional equation obtained by constructing
quadrant walks step by step:
\beq\label{eq-Q0}
K(x,y)Q(x,y)=1-t\by Q_+(x)-t\bx Q_+(y)
\eeq
where 
$$
Q_+(x)=Q(x,0)=\sum_{i\ge 0, n\ge 0} q_{i,0}(n) x^i t^n \ \ \in
\qs[x][[t]].
$$
Equivalently,
\beq\label{eq-Q}
xyK(x,y)Q(x,y)=xy-txQ_+(x)-tyQ_+(y).
\eeq

The kernel $K(x,y)$ is invariant by the transformations $x\mapsto \bx$
and $y\mapsto \by$.   Hence we also have:
\begin{align*}
  \bx yK(x,y)Q(\bx,y)&=\bx y-t\bx Q_+(\bx)-tyQ_+(y),\\
\bx \by K(x,y)Q(\bx,\by)&=\bx \by-t\bx Q_+(\bx)-t\by Q_+(\by),\\
x \by K(x,y)Q(x,\by)&=x \by-txQ_+(x)-t\by Q_+(\by).
\end{align*}
The  \emm orbit equation, is the alternating sum of the last
four equations:
\begin{align}
K(x,y) \left(xyQ(x,y)-\bx y Q(\bx,y)+\bx \by Q(\bx, \by)-x\by
  Q(x,\by)\right)&=xy-\bx y +\bx \by -x\by \nonumber\\
&=(x-\bx)(y-\by).\label{orbit-Q}
\end{align}
Observe that the right-hand side is now explicit. We call it the \emm
orbit sum, of this quadrant model.
The above equation can be rewritten as
\beq\label{Q-lin}
xyQ(x,y)-\bx y Q(\bx,y)+\bx \by Q(\bx, \by)-x\by Q(x,\by)=
\frac{(x-\bx)(y-\by)}{1-t(x+\bx+y+ \by)}.
\eeq
Extracting the positive part in $x$ and $y$ gives the following
classical result~\cite{guy-bijections,BoMi10}.
\begin{prop}
  The series $xyQ(x,y)$ is the positive part (in $x$ and $y$) of the
  rational function
$$
\frac{(x-\bx)(y-\by)}{1-t(x+\bx+y+ \by)}.
$$
It is thus D-finite.  For $i,j,m\ge 0$, the number of quadrant walks
of length $n=i+j+2m$ starting at $(0,0)$ and 
ending at $(i,j)$ is:
$$
\frac{(i+1)(j+1)n!(n+2)!}{m!(m+i+j+2)!(m+i+1)!(m+j+1)!}.
$$
\end{prop}
\subsection{Reduction to an equation with orbit sum zero}
We now return to walks avoiding the negative quadrant. We first apply
to the functional equation~\eqref{eq-T}, written as
$$
xyK(x,y) \C(x,y)= xy- tx \Cn(\bx) -ty \Cn(\by),
$$
 the  treatment that we have just applied to the quadrant equation~\eqref{eq-Q}.
 The orbit equation reads as above
\beq\label{orbit-C}
K(x,y) \left( xy\C(x,y)-\bx y \C(\bx,y)+\bx \by \C(\bx, \by)-x\by
  \C(x,\by)\right) =(x-\bx)(y-\by).
\eeq
We have just seen that $Q(x,y)$ also satisfies this equation
(see~\eqref{orbit-Q}). The same holds for $-\bx^2 Q(\bx,y)$, for $-\by^2 Q(x,\by)$
 and for $\bx^2\by^2Q(\bx,\by)$ (by~\eqref{orbit-Q} again). 
Let $a$ 
be a real number and write
\beq\label{A-def}
\C(x,y)=A(x,y) +(1-2a)Q(x,y) -a \bx^2 Q(\bx,y) -a \by^2 Q(x,\by),
\eeq
where $A(x,y)$ is a new series. (We do not involve the fourth solution $\bx^2\by^2Q(\bx,\by)$,
 since $\C(x,y)$ contains no monomial that is negative in $x$ and in
 $y$). 
The above observations imply that the orbit sum associated with
$A$ vanishes:
\beq\label{orbit-A}
xyA(x,y)-\bx y A(\bx,y)+\bx \by A(\bx, \by)-x\by A(x,\by)=0.
\eeq
Moreover, we can compute from the equations~\eqref{eq-T} and~\eqref{eq-Q0}
satisfied by $\C$ and $Q$ a functional equation for $A$.
We first note that 
$$
\Cn(\bx)=A_-(\bx) -a \bx^2 Q_+(\bx),
$$
where, as above,
$$
A_-(\bx)= \sum_{i<0, n\ge 0} a_{i,0}(n) x^i t^n \ \ \in \bx \qs[\bx][[t]].
$$
Then we derive  from~\eqref{eq-T} and~\eqref{eq-Q0} that:
$$
K(x,y)A(x,y)= a(2+\bx^2+\by^2) -t\by A_-(\bx) -t \bx A_-(\by)
+t\by (1-3a) Q_+(x)
+t\bx(1-3a) Q_+(y).
$$
This suggests to choose $a=1/3$, so that
\beq  \label{eq-A}
K(x,y)A(x,y)= (2+\bx^2+\by^2)/3-t\by A_-(\bx) -t \bx A_-(\by) .
\eeq
Note that the equations~\eqref{eq-T} and~\eqref{eq-A} satisfied by $\C$ and $A$ only differ by their
constant term on the right-hand side: it is simply $1$ for $\C$, but
$(2+\bx^2+\by^2)/3$ for $A$. This results into a zero orbit sum for~$A$. Observe also that~\eqref{eq-A} characterizes $A(x,y)$ uniquely as
a \fps\ in $t$. The series $3A(x,y)$ counts walks in $\cC$ starting
from $(0,0)$, $(-2,0)$ or $(0,-2)$, but those starting at $(0,0)$ get
weight 2.

We will show that  $A(x,y)$ is algebraic (hence the notation
$A$), as claimed by Theorem~\ref{thm:square}. We recall that all quadrant models with small steps and orbit sum
zero have an algebraic
\gf\ too~\cite{BoKa10,BoKuRa13,mbm-gessel,bousquet-versailles,Bous05,BoMi10,gessel-proba}. One
main difference with these quadrant models
is that $A(x,y)$ involves positive \emm and negative, powers of $x$ and~$y$. The next
subsection takes care of this difference.

\subsection{Reduction to a quadrant-like problem}
\label{sec:quadrant-red}
We now separate in $A(x,y)$ the contribution of the three quadrants:
\beq\label{A-expr}
A(x,y)=P(x,y)+ \bx M(\bx,y)+\by M( \by,x),
\eeq
where $P(x,y) \in \qs[x,y][[t]]$ and $M(\bx,y) \in
\qs[\bx,y][[t]]$. Note that we have exploited the obvious $x/y$
symmetry, and that this identity defines $P$ and $M$ uniquely in terms
of $A$. The
letter $P$ stands for \emm positive,, and the letter $M$ for \emm
mixed,. 
Extracting the positive part in $x$ and $y$ from the orbit
equation~\eqref{orbit-A}  gives
\beq\label{PM-bis}
xyP(x,y)=y\left( M(x,y)-M(0,y)\right) + x\left(
M(y,x)-M(0,x)\right),
\eeq
which is equivalent to~\eqref{PM}. Hence it suffices to determine $M$. 
We can  write the above series $A$ in terms of $M$ only:
$$
A(x,y)=\bx\left( M(x,y)-M(0,y)\right) + \by \left(
M(y,x)-M(0,x)\right)+ \bx M(\bx,y)+\by M( \by,x).
$$
Note that
$
A_-(\bx)= \bx M(\bx,0).
$
Plugging this in the
functional equation~\eqref{eq-A} gives:
\begin{multline*}
  K(x,y) \left( \bx M(x,y)-\bx M(0,y)+\by M(y,x)-\by
M(0,x) +\bx M(\bx,y) + \by M(\by, x)\right)
\\ =(2+\bx^2+\by^2)/3 -t\bx\by M(\bx,0)-t\bx\by
M(\by,0).
\end{multline*}
Let us  extract  the negative part in $x$ of this equation. We obtain:
\begin{multline}\label{M-neg}
  -t\bx \big( M_x(0,y)  + \by M(y,0)-\by M(0,0) \big)
+K(x,y) \bx M(\bx,y)+ t M(0,y) -t\bx \by M(\by,0)\\= \bx^2/3 -t\bx\by M(\bx,0)-t\bx\by M(\by,0),
\end{multline}
with $M_x=\partial M/\partial x$. Observe that the term $t\bx \by
M(\by,0)$ occurs on both sides and thus cancels.
Extracting from this the coefficient of $\bx$ gives:
$$
-t\left( M_x(0,y)+ \by M(y,0)-\by M(0,0) \right)
+(1-t(y+\by)) M(0,y)-tM_x(0,y)= -t\by M(0,0).
$$
From this, we obtain an expression of $ M_x(0,y)$ in terms of
$M(y,0)$  and $M(0,y)$. By plugging it in~\eqref{M-neg},  we obtain
$$
K(x,y) \left( 2M(\bx,y)-M(0,y)\right) = 2\bx/3 -2t \by M(\bx,0) + t(\bx-x)
M(0,y)  +t\by M(y,0).
$$
 Replacing $x$ by $\bx$ gives the  functional
equation~\eqref{func-M} for $M(x,y)$, which we repeat here for convenience:
\beq\label{eqcat-M}
K(x,y) \left( 2M(x,y)-M(0,y)\right) = 2x/3 -2t \by M(x,0) + t(x-\bx)
M(0,y)  +t\by M(y,0).
\eeq
We recall that $M(x,y)$ is a series in $t$ with polynomial coefficients
in $x$ and $y$. The equation it satisfies is reminiscent of the quadrant
equation~\eqref{eq-Q0}. However, its right-hand side involves the series
$M(y,0)$ in addition to the two standard specializations  $M(x,0)$ and
$M(0,y)$. 
Still, several ingredients in the rest of the solution are borrowed
from  former solutions of quadrant models~\cite{bousquet-versailles,mbm-gessel,BoMi10}.

\subsection{Cancelling the kernel: an equation between
  $\boldsymbol{M(x,0)}$, $\boldsymbol{M(0,x)}$ and
  $\boldsymbol{M(0,\bx)}$}
\label{sec:cancel}
As a polynomial in $y$, the kernel $K(x,y)$ has two roots. Only one of
them is a power series in $t$ (with coefficients in $\qs[x,\bx]$). We
denote it by $Y\equiv Y(x)$:
\beq\label{Y-expr}
Y=\displaystyle \frac{1-t(x+\bx) -\sqrt
{ (1-t(x+\bx))^2-4t^2}
}{2t}=
t+(x+\bx)t^2+O(t^3) .
\eeq
The other root is $1/Y$, and its expansion in $t$ involves a term $1/t$. 

Specializing $y$ to $Y$ in~\eqref{eqcat-M} gives a relation between
the three series on the right-hand side:
\beq\label{kernel1}
2x/3 -2t  M(x,0)/Y + t(x-\bx) M(0,Y)  +tM(Y,0)/Y=0.
\eeq
Since $Y$ is symmetric in $x$ and $\bx$, we also have, upon replacing
$x$ by $\bx$:
\beq\label{kernel2}
2\bx/3 -2t  M(\bx,0)/Y + t(\bx-x)
M(0,Y)  +tM(Y,0)/Y=0.
\eeq
Since the kernel is symmetric in $x$ and $y$, we can
rewrite~\eqref{eqcat-M} as
$$
K(x,y) \left( 2M(y,x)-M(0,x)\right) = 2y/3 -2t \bx M(y,0) + t(y-\by)
M(0,x)  +t\bx M(x,0).
$$
Setting $y=Y$ in this equation gives
\beq\label{kernel3}
2Y/3 -2t \bx M(Y,0) + t(Y-1/Y)M(0,x)  +t\bx M(x,0)=0.
\eeq
Finally, using once more  $Y(x)=Y(\bx)$ gives a fourth
equation:
\beq\label{kernel4}
2Y/3 -2t x M(Y,0) + t(Y-1/Y)M(0,\bx)  +txM(\bx,0)=0.
\eeq
We have thus obtained four
equations, namely~(\ref{kernel1}--\ref{kernel4}), relating the six
series $M(x,0)$, $M(\bx,0)$, $M(Y,0)$, $M(0,x)$, $M(0, \bx)$ and $M(0,
Y)$. By eliminating $M(0,Y)$ and $M(Y,0)$, we obtain 
two equations between the remaining series, which only differ by the transformation $x\mapsto
\bx$. Eliminating  $M(\bx,0)$ between them gives:
\beq\label{catM}
(Y-1/Y) \left( xM(0,x)-2\bx M(0,\bx)\right)- 2\bx Y/t +3 M(x,0)=0.
\eeq
Our next step is to eliminate $M(x,0)$.

\subsection{An equation between $\boldsymbol{M(0,x)}$ and
  $\boldsymbol{M(0,\bx)}$}
\label{sec:M0x-M0bx}

Let us denote the discriminant occurring in the expression~\eqref{Y-expr} of $Y$ by
$$
\Delta(x):= { (1-t(x+\bx))^2-4t^2}.
$$
Denote also
\beq\label{RS-def}
R(x)= tM(x,0) \qquad \hbox{and} \qquad S(x)=txM(0,x).
\eeq
Then the above equation~\eqref{catM} reads
\beq\label{eqRS}
\sqrt{\Delta(x)}\left( S(x)-2S(\bx)-\bx\right)
+\bx -t-t\bx^2=3tR(x).
\eeq
Consequently,
$$
{\Delta(x)}\big( S(x)-2S(\bx)-\bx\big)^2=\big(3tR(x)
-\bx+t+t\bx^2
\big) ^2.
$$
Recall that $S(x)$ is a series in $t$ with coefficients in $x\qs[x]$. Extracting from the
above identity the negative part in $x$ gives (after dividing by $4$):
$$
  \Delta(x) \left ( S(\bx)^2+\bx S(\bx) \right) 
- [x^{<}]\left(  \Delta(x) S(x) S(\bx)\right) =
F_0 +\bx F_1 + \bx^2 F_2,
$$
where $F_0, F_1$ and $F_2$ are series in $t$ than can be expressed in
terms of  $S(x)$ and $R(x)$:
\beq\label{F2-expr}
F_0 = t^2\Sp(1+\Sp), \qquad 
F_1=\frac{ t}2 \left( tS_2 +3  tR_1-5\Sp\right), 
\qquad F_2= t^2(1+2\Sp),
\eeq
where we have used the notation~\eqref{Gi-def} 
(note that $R_0=S_1$ by definition of $R$ and $S$).
We have thus obtained an expression for the negative part of
the series $\Delta(x)S(x)S(\bx)$:
$$
   [x^{<}]\left(  \Delta(x) S(x) S(\bx)\right) 
=
\Delta(x) \left ( S(\bx)^2+\bx S(\bx) \right)-F_0 -\bx F_1 - \bx^2
F_2.
$$
Let us denote 
\beq\label{P0-def}
P_0:=[x^0] \left(\Delta(x)S(x)S(\bx)\right).
\eeq
  By symmetry, we can now express $\Delta(x)S(x)S(\bx)$ as follows:
$$
 \Delta(x) S(x) S(\bx)=P_0+
\Delta(x) \left ( S(x)^2+S(\bx)^2+xS(x)+\bx S(\bx) \right) 
-2F_0-(x+\bx)F_1-(x^2+\bx^2) F_2.
$$
Equivalently:
\beq
\label{eqcat-S}
  \Delta(x) \big(S(x)^2+S(\bx)^2-S(x)S(\bx) +xS(x)+\bx S(\bx)\big)=
2F_0-P_0+(x+\bx)F_1+(x^2+\bx^2) F_2.
\eeq
This is the promised equation relating $S(x)$ and $S(\bx)$, or
equivalently, $M(0,x)$ and $M(0, \bx)$. In the next step, we will get
rid of $M(0,\bx)$.

\medskip
\noindent{\bf Remark.} One may wonder whether the right-hand side of~\eqref{eqcat-S}, as
a polynomial in $x$, is divisible by $\Delta(x)$, or at least by one
of its two factors $(1-t(x+\bx+2))$ and $(1-t(x+\bx-2))$. This would
 simplify our calculations, but one can easily check that this is not
 the case.

\subsection{An equation for  $\boldsymbol{M(0,x)}$ only}
\label{sec:M0x}
The product $S(x)S(\bx)$ occurring in the above equation~\eqref{eqcat-S} makes it
difficult  to extract the positive part. To eliminate this cross term,
we multiply the equation by $S(x)+S(\bx)+x+\bx$. This ``trick'' (based
on the identity $(a^2+b^2-ab)(a+b)=a^3+b^3$) was
already used in the solution~\cite{mbm-gessel} of Gessel's quadrant model. This
gives:
\begin{multline*}
  \Delta(x) \big( S(x)^3+ S(\bx)^3+ (2x+\bx) S(x)^2 + (2\bx+x)
  S(\bx)^2 + x(x+\bx) S(x) +\bx(x+\bx)S(\bx) \big) 
   \\=
\left(2F_0-P_0+(x+\bx)F_1+(x^2+\bx^2) F_2\right)
 \left(S(x)+S(\bx)+x+\bx\right).
\end{multline*}
Now we can extract the non-negative part in $x$:
\begin{multline*}
  \Delta(x) \big( S(x)^3+ (2x+\bx) S(x)^2 
+ x(x+\bx) S(x)\big) = t^2(x-\bx)(1+\Sp)^2+(1+\Sp)(F_1+2t\Sp)\\
+\left(2F_0-P_0+(x+\bx)F_1+(x^2+\bx^2) F_2\right)\left(S(x)+x\right).
\end{multline*}
(We have used the expression~\eqref{F2-expr} of
$F_2$ to simplify the right-hand side).
Extracting the
constant term in $x$ gives $F_1=-2t\Sp $, 
so that the equation satisfied by $S(x)=txM(0,x)$ is:
\begin{multline}\label{eqcatS}
  \Delta(x) \left( S(x)^3+ (2x+\bx) S(x)^2 + x(x+\bx) S(x)\right) = 
t^2(x-\bx)(1+\Sp)^2\\
+\left(2t^2\Sp ^2 +2t\left(tx^2+t\bx^2 -x-\bx+t \right) \Sp   -P_0
  +t^2(x^2+\bx^2)\right)\left( S(x)+x\right).
\end{multline}
In addition to $S(x)$, it involves two  series depending on $t$
only, namely $\Sp$ and $P_0$. Still, we will see that this equation, combined with the
fact that $S(x)$ has polynomial coefficients in $x$ and the values of
the first few of these coefficients, determines uniquely $P_0$ and $S(x)$ 
 (and consequently $S_1$).

\subsection{The generalized quadratic method: algebraicity of $\boldsymbol{A(x,y)}$}
\label{sec:quad}
\paragraph{\bf General principle.} We have described in~\cite{mbm-jehanne} how to study  equations of the form
\beq\label{quad1}
\Pol(S(x), A_1, \ldots, A_k,t,x)=0,
\eeq
where $\Pol(x_0, x_1, \ldots, x_k, t,x)$ is a polynomial with
 complex coefficients, $S(x)$ is a \fps\ in $t$ with coefficients in
$\qs[x]$, and $A_1, \ldots, A_k$ are $k$ auxiliary series depending
on  $t$ only, under the assumption that these $k+1$ series are uniquely determined by~\eqref{quad1}. (In the above example~\eqref{eqcatS}, $\Pol$ is a \emm Laurent, polynomial in
 $x$, but this makes no difference.)
The strategy of~\cite{mbm-jehanne} instructs us to look for power series
$X \in \cs[[t]]$ satisfying
\begin{equation}
\frac{\partial \Pol }{\partial x_0}(S(X), A_1, \ldots, A_k,t,X)=0.
\label{quad2} \end{equation}
Indeed, by differentiating~\eqref{quad1} with respect to $x$, we see
that any such series also satisfies
\begin{equation}  
 \frac{\partial \Pol }{\partial x}(S(X), A_1, \ldots, A_k,t,X)=0,
\label{quad3} \end{equation}          
and we thus obtain  three polynomial equations, namely
Eq.~\eqref{quad1} 
written for $x=X$, Eqs.~\eqref{quad2} and~\eqref{quad3}, that relate
the $(k+2)$ unknown series $S(X)$, $A_1, \ldots, A_k$ and~$X$.
 If we
can prove the existence of $k$ distinct series $X_1, \ldots, X_k$
satisfying~\eqref{quad2}, we will 
have $3k$ equations between the $3k$ unknown series $S(X_1), \ldots ,
S(X_k)$, $A_1, \ldots, A_k$, $X_1, \ldots, X_k$. If there is no redundancy
in this system, we will have proved 
that each of the $3k$ unknown series is algebraic over
$\cs(t)$. 

\medskip
\paragraph{\bf Identifying the series $\boldsymbol{X_i}$.} We apply this strategy to~\eqref{eqcatS}, with $A_1=\Sp $ and
$A_2=P_0$. The polynomial $\Pol(x_0, x_1, x_2, t,x)$ is
\begin{multline}\label{Pol-def}
  \Delta(x) \left( x_0^3+ (2x+\bx) x_0^2 + x(x+\bx) x_0\right) -t^2(x-\bx)(1+x_1)^2 \\
-\left(2t^2x_1^2 +2t\left(tx^2+t\bx^2 -x-\bx+t \right) x_1  -x_2
  +t^2(x^2+\bx^2)\right) \left(x_0+x\right).
\end{multline}
Equation~\eqref{quad2} reads:
\begin{multline}\label{quad2G}
 \Delta(X)\big( 3S(X)^2+2(2X+1/X)S(X)+X(X+1/X)\big) 
\\= 2t^2\Sp ^2 +2t\left(tX^2+t/X^2 -X-1/X+t \right) \Sp   -P_0
  +t^2(X^2+1/X^2).
\end{multline}
Recall the definitions~\eqref{RS-def} and~\eqref{P0-def} of $S$ and
$P_0$.   Once multiplied by $X$,  the above equation has the following
form: 
$$
X(1+X^2)=-X P_0+ t \Pol_1\left(\frac{M(0,X)-M(0,0)}X, M(0,0), t,X\right)
$$
for some polynomial $\Pol_1$.
This already shows that there exist two solutions $X_1$ and $X_2$ in
$\cs[[t]]$ having constant terms $i$ and $-i$ respectively: since
$P_0$ is a multiple of $t$, the above equation allows one to compute
their coefficients inductively, assuming $M(0,x)$ is known. But we will also use a third solution,
which has constant term $0$. To show its existence, let us write $X=t
\hat X$. Then the above equation, once divided by $t$, reads
$$
\hat X=2- \hat X P_0+ t \Pol_2\left(\frac{M(0,t\hat X)-M(0,0)}{t\hat
    X}, M(0,0), t,\hat X\right),
$$
which gives  a third solution $X_0$, of the form $2t+O(t^2)$.
Using the first few coefficients of $M(0,x)$, we  obtain: 
\begin{align}
  X_0&=2t+8t^3+64t^5+640t^7+7168 t^9+O(t^{11}),\label{init-X0}
\\
X_{1,2}&=\pm i + 2t^3+16 t^5 \mp 2i t^6 +156 t^7 +O(t^8).\label{init-X12}
\end{align}
In particular, these three series are non-zero, and  are thus the three solutions of~\eqref{quad2G}. 

We note that 
$$
X_0=\frac 1{2t} \left( (4t^2) + (4t^2)^2 +
2(4t^2)^3 +5(4t^2)^4+14(4t^2)^5 +\cdots\right),
$$
seems to be related to Catalan numbers. Due to the special
 form of our polynomial $\Pol$ (given by~\eqref{Pol-def}), it is in
 fact simple to prove that
 \beq\label{X0-sol}
X_0= \frac{1-\sqrt{1-16t^2}}{4t}.
\eeq
 Indeed, we observe that 
\begin{multline*}\label{simple}
(x^2+1) \Pol -x(x^2-1)\frac{\partial \Pol}{\partial x} -
(2x+x_0+x^2x_0)  \frac{\partial 
  \Pol}{\partial x_0}=
 -2x(1-2t(x+\bx))\\ \times (x+x_0)
 \left( x+\bx +t(x-\bx)^2x_1+x_0(x_0+x+\bx)\left(x+\bx-t(x-\bx)^2\right)\right).
\end{multline*}
 Since the three series $X\equiv X_i$ cancel $\Pol$ and its partial
 derivatives, each of them  must satisfy 
\beq\label{fact-1}
1-2t\left( X+1/X\right)=0,
\eeq
or
$$
X+S(X)=0,
$$
or
\beq\label{fact-3}
X+1/X +t(X-1/X)^2 \Sp +S(X)(S(X)+X+1/X)\left(X+1/X-t(X-1/X)^2\right)=0.
 \eeq
Using the initial values~(\ref{init-X0}-\ref{init-X12}), we conclude that $X_0$
satisfies~\eqref{fact-1} (from which the expression~\eqref{X0-sol}
follows), and that $X_1$ 
and $X_2$ satisfy~\eqref{fact-3}.

\medskip
\paragraph{\bf Elimination.}
Eliminating $S(X)$ and $X$ between $\Pol(S(X),\Sp , P_0,t,X)$,
$\Pol_{x_0} (S(X), \Sp , P_0, t, X)$ and \eqref{fact-1} gives a first
polynomial equation between $\Sp $ and $P_0$. A second equation is
obtained by eliminating $S(X)$ and $X$ between $\Pol(S(X),\Sp ,
P_0,t,X)$, 
$\Pol_{x_0} (S(X),\Sp , P_0,t,X)$ and~\eqref{fact-3}. (When several
factors occur, one determines the correct one using the first
coefficients of $S(x)$, $\Sp $, $P_0$ and the $X_i$'s.) A further
elimination (first of $P_0$, then of $\Sp $) between these two
equations gives polynomial equations for each of these two
series. Both are found to be of degree 4 over $\qs(t)$:
\begin{multline}\label{S1-eq-sq}
  19683{t}^{6}\Sp ^{4}+2187{t}^{4} \left( 20{t}^{2}-1
 \right) \Sp ^{3}+81{t}^{2} \left( 11{t}^{2}-1 \right) 
 \left( 38{t}^{2}-1 \right) \Sp  ^{2}\\
+ \left( 92{t}^{2}-1
 \right)  \left( 11 t^{2}-1 \right) ^{2}\Sp +{t}^{2} \left( 
1331{t}^{4}-107{t}^{2}+1 \right) =0,
\end{multline}
and
\begin{multline*}
  387420489{t}^{6}P_0^{4}+3188646{t}^{4} \left( 284{t}^{4}
-113{t}^{2}-1 \right) P_0^{3}\\+8748{t}^{2} \left( 31570{t
}^{8}-96755{t}^{6}+7251{t}^{4}+{t}^{2}+1 \right) P_0^{2}\\+
 \left( 29962144{t}^{12}-441273288{t}^{10}+87261432{t}^{8}-
4754122{t}^{6}+64860{t}^{4}-687{t}^{2}-8 \right) P_0\\+{t}^{
4} \left( 1102736{t}^{10}-53770928{t}^{8}+4286896{t}^{6}-58740
{t}^{4}+751{t}^{2}+8 \right) =0.
\end{multline*}
From this and~\eqref{eqcatS}, we derive that $S(x)=txM(0,x)$ is algebraic
over $\qs(t,x)$. Then,~\eqref{catM} implies that
the same holds for $M(x,0)$. Finally, the algebraicity of
$M(x,y)$ follows from~\eqref{eqcat-M}, and that of $P(x,y)$
from~\eqref{PM-bis}. We have thus proved the algebraicity of the series
$A(x,y)$ given by~\eqref{A-expr}, which, by definition, is
$$
\C(x,y)-\frac 1 3 Q(x,y) +\frac 1 3 \bx^2 Q(\bx,y) +\frac 1 3 \by^2
Q(x,\by) .
$$

\subsection{Rational parametrizations and degrees}
\label{sec:param}
The equations obtained above for $\Sp $ and $P_0$  can be
parametrized by introducing the unique series 
$T\in \qs[[t]]$, with constant term $1$, satisfying~\eqref{T-def}.
Indeed, both equations factor when replacing $t^2$ by
$(T-1)(T+3)^3/(256 T^3)$,
and extracting the correct factor gives:
\beq\label{Sp0sol}
\Sp                    =
 \frac{(T-1)(11+6T-T^2)}{(T+3)^3}
\eeq
and
$$
P_0
=
\frac{(T-1)^2(41+331T+106T^2+38T^3-3T^4-T^5)}{128T^3(T+3)^3}.
$$
We recall that parametrizations of algebraic curves (of genus 0) can be computed
using  the  Maple command {\tt parametrization}.
 
We now plug these expressions in the equation~\eqref{eqcatS} defining
$S(x)$. This gives a cubic equation for $S(x)$ over
$\qs(t,x,T)$. Eliminating $T$ gives an irreducible polynomial of
degree 12 
in $S(x)$ over $\qs(t,x)$: since $S(x)=txM(0,x)$, we conclude that
$M(0,x)$ has degree 12 as well. 

We now return to the cubic equation satisfied by $S(x)$ over
$\qs(t,x,T)$, and  replace $x$ by $xt$. Due to the structure of the
square lattice,  $S(xt)/t=xtM(0,xt)$  is an
even function of~$t$: this allows us to replace $t^2$ by its
rational expression in terms of $T$, and  gives a cubic equation
for $S(xt)/t$ over $\qs(x,T)$. Then, introducing the parametrization~\eqref{param-x} of $x$ factors 
this cubic equation into a linear factor and a quadratic one. The one that
vanishes is found to be the linear one. This gives a rational
expression of $S(xt)/t$ in terms of $T$ and $U$, which is equivalent
to the expression~\eqref{M0xt} of $tM(0,xt)$.

We now want to express $M(x,0)$, or equivalently the series
$R(x)=tM(x,0)$, using~\eqref{eqRS}. Using the cubic equation (over $\qs(t,x,T)$) found for
$S(x)$, we first find that the term $D(x):=(S(x)-2S(\bx)-\bx)$ has degree 6
over $\qs(t,x,T)$, but is in fact bicubic. Thus~\eqref{eqRS} gives
automatically an equation of degree 6 for $R(x)$ over
$\qs(t,x,T)$. Eliminating $T$ shows that $R(x)$ has degree 24 over
$\qs(t,x)$. Since $R(x)=tM(x,0)$, the same holds for $M(x,0)$.

Finally,   we replace $x$ by $xt$ in the equation of degree 6
satisfied by $R(x)$ over $\qs(t,x,T)$  (the series $R(xt)$ is an even
function of $t$).  Then we parametrize $t^2$ and $T$ by $Z=\sqrt T$, and find that $R(xt)$ is cubic over $\qs(x,Z)$. We finally parametrize $x$ by $U$ (as in~\eqref{param-x}): the equation
factors into a linear term and a quadratic one. The one that
cancels turns out to be linear, and this gives for $R(xt)$ a
rational expression in $Z$ and $U$ which is equivalent to our
expression~\eqref{Mxt0} of $tM(xt,0)$.

\medskip
It remains to show that $M(x,y)$ and $P(x,y)$ have degree 72 over
$\qs(x,y,t)$. It suffices to prove this for $tM(xt,yt)$ and $P(xt,yt)$,
which are even series in $t$.

First, it follows from~\eqref{eqcat-M} and the rational expressions of
$tM(0,xt)$ and $tM(xt,0)$ that $tM(xt,yt)$ belongs to
$\qs(Z,U,\tilde U)$, where $\tilde U$ is the counterpart of $U$ for
the variable $y$ instead of $x$. Since $Z$ has degree 8 over $\qs(t)$,
and $U$ has degree 3 over
$\qs(x,Z)$, it follows that $tM(xt,yt)$ has degree at most 72 over
$\qs(t,x,y)$. Computing its minimal equation over $\qs(x,y,t)$ (in
practice, for $x=2$ and $y=-2$ for instance, to avoid extremely heavy
computations) shows that this bound is tight.

We proceed similarly with $P(x,y)$, expressed in terms of $M(x,y)$ and
its specializations thanks to~\eqref{PM-bis}.

We have now completed the proof of Theorem~\ref{thm:square}.

\subsection{Walks ending at a prescribed position}
\label{sec:prescribed}
Let us now prove Corollary~\ref{cor:coeffs}. We want to show that for
$i, j \ge 0$, the
coefficients of $x^iy^j$ in $M(x,y)$ and $P(x,y)$, which we denote by
$M_{i,j}$ and $P_{i,j}$ respectively,  belong to $t^{i+j+1} \qs(Z)$
(resp. to $t^{i+j} \qs(Z)$). First, the
connection~\eqref{PM} between $P$ and $M$ gives
$$
P_{i,j}=M_{i+1,j}+M_{j+1,i},
$$
and shows that it suffices to prove
the property for the series $M_{i,j}$. Then, extracting from the functional equation~\eqref{func-M} satisfied
by $M$ the coefficient of $x^i y^j$ 
shows that for $i, j\ge 0$, the series $M_{i,j+1}$ can be expressed as a
linear combination of series $M_{k,\ell}$ such that $\ell \le j$
and/or $k=0$, with coefficients in $\qs[t,1/t]$. 
Hence it suffices to prove our
results for the series $M_{i,0}$ and $M_{0,j}$.
Equivalently, by definition~\eqref{RS-def} of the series $R$ and $S$, it suffices to
prove that the coefficient of $x^i$ in $\tilde
R(x):= R(xt)=tM(xt,0)$ and $\tilde S(x):= S(xt)/(tx)=tM(0,xt)$   (which are
both even functions of $t$) belong to $\qs(Z)$. 

Recall from the previous subsection that $\tilde R$ satisfies a cubic
equation over 
$\qs(Z,x)$. This equation reads
\beq\label{eqR}
Z^{24} (Z^2+3)^3 \tilde R = Z^{24}(Z^2-1)(11+6Z^2-Z^4) + x \Pol(x,Z,
\tilde R)
\eeq
for some polynomial $\Pol$.
This gives 
$$
[x^0] \tilde R = [x^0] R(x) = t M_{0,0} =
\frac{(Z^2-1)(11+6Z^2-Z^4)}{(Z^2+3)^3}
$$
(which we already obtained in~\eqref{Sp0sol}), and shows, by induction on
$i\ge 0$, that the coefficient of $x^i$ in $\tilde R(x)$ is a rational
function of $Z$.

For the series $\tilde S(x)=tM(0,xt)$, we have to go one step further in the
application of Newton's polygon method. We start from the cubic
equation satisfied by this series over $\qs(T,x)$, which we write as
an equation over $\qs(Z,x)$ with $T=Z^2$. Writing $\tS(x)=\tS_0+x\tS_1+
x^2 \hat S(xà)$, we first determine (by setting $x=0$ in the equation) the
value of $\tS_0$ (once again equivalent to~\eqref{Sp0sol}). Then there are
two possible values for $\tS_1$; after checking the first
coefficients, we conclude that the correct one is
$$
\tS_1= t^2M_{0,1}=\frac{4(Z-1)^2(Z^2+1)^2(1+2Z-Z^2)}{Z^3(3+Z^2)^3}.
$$
 For the remaining series $\hat S(x)$, we find an equation
over $\qs(Z,x)$  
which, as~\eqref{eqR}, has degree~1 in $\hat S$ when $x=0$. One can
then  compute recursively the coefficient of $x^i$ in $\hat S(x)$,
which belongs to $\qs(Z)$.

Finally, the nature of the series $C_ {i,j}$ follows from the fact that
$Q_{i,j}$ is D-finite but transcendental for $i, j \ge 0$. Indeed, the
asymptotic behaviour its $n$th coefficient, in
$4^nn^{-3}$, is not compatible with algebraicity~\cite{flajolet-context-free}.

\section{The diagonal square lattice}
\label{sec:diag}
We now adapt the calculations of Section~\ref{sec:square} to
walks on the \emm diagonal, square
lattice 
(Figure~\ref{fig:square-diag}, right). That is,
walks now take steps $(\pm 1, \pm 1)$. One difference from the square
lattice case is an extra term in the basic functional equation,
corresponding to the forbidden move from $(0,0)$ to
$(-1,-1)$. Otherwise the argument is very similar, and we  give
much fewer details.

We adopt the same notation as before. In particular, $C(x,y)$ 
denotes the \gf\ of walks on the diagonal square lattice, starting
from $(0,0)$ and  confined to 
$\cC$. 
As in the square lattice case, the expression of $C(x,y)$  involves the \gf \
$Q(x,y)$ of walks confined to the first quadrant, which is now~\cite{BoMi10}:
\beq\label{Q-diag}
Q(x,y)= \sum_{i,j,n\ge 0}
\frac{(i+1)(j+1)}{\left(1+\frac{n+i}2\right)\left(1+\frac{n+j}2\right)}
  \binom{n}{ \frac{n+i}2} \binom{n}{ \frac{n+j}2}\, x^iy^jt^n,
\eeq
where the sum is restricted to values of $i$, $j$ and $n$ having the same parity.

\begin{theorem}\label{thm:diag}
 The \gf\ of walks with steps $(\pm 1, \pm 1)$, starting at $(0,0)$,
 confined to $\mathcal \C$  
and ending in the first quadrant
(resp. at a negative abscissa) is
$$
\frac 1 3\, Q(x,y) + P(x,y), \qquad\qquad
\left( \hbox{resp.} \quad -\frac 1 3\, \bx^2 Q(\bx,y)+  \bx M(\bx,y)
\right),
$$
where $M(x,y)$ and $P(x,y)$ are algebraic series of degree $72$ over $\qs(x,y,t)$.

More precisely,  $P$ can be expressed in terms of $M$ by:
\beq\label{PM-diag}
P(x,y)=\bx\left( M(x,y)-M(0,y)\right) + \by \left(
M(y,x)-M(0,x)\right),
\eeq
$M$ satisfies the functional equation
\begin{multline}\label{func-M-diag}
  \left(1-t(x+\bx)(y+\by)\right)\left( 2M(x,y)-M(0,y)\right) =
2x/3-2t\by(x+\bx)M(x,0)\\+t(x-\bx)(y+\by)M(0,y)+t(1+\by^2)M(y,0)-t\by M_{1,0}
\end{multline}
where $M_{1,0}$ is the coefficient of $x^1y^0$ in $M(x,y)$, and the specializations $M(x,0)$ and $M(0,y)$ have respective degrees
$24$ and $12$  over $\qs(t,x)$ and $\qs(t,y)$. 

Moreover, these algebraic  series admit
rational parametrizations. Let us define the series $T$ and $Z$ as in Theorem~{\rm\ref{thm:square}}, and let
 $V$ be the only series in $t$, with constant term $0$,
satisfying 
\beq\label{V-def}
1-T+3V+VT= xV^2(3+V+T-VT).
\eeq
Then
  the series $\sqrt x M(\sqrt x,0)$ and $t\sqrt x M(0,\sqrt
x)$ (both even series in $t$ with polynomial coefficients in $x$)
admit rational expressions in terms of $Z$ and $V$, given in Appendix~\ref{sec:param-diag}.
\end{theorem}
By combining the above results and \emm singularity analysis, of
algebraic (and D-finite) series~\cite{flajolet-sedgewick}, one can
derive from the above theorem that the number
of $n$-step walks on the diagonal square lattice confined to $\cC$ is
$$
[t^n] C(1,1)\sim \frac{2^3 \sqrt 3}{3^2\, \Gamma(2/3)}
\frac{4^n}{n^{1/3}}.
$$ 

\begin{cor}[\bf Walks ending at a prescribed position]\label{cor:coeffs-diag}
 Let $T$ be the unique series in $t$ with constant term $1$
  satisfying~\eqref{T-def}, and let  $Z=\sqrt T $.

 For $j\ge 0$, the series $C_{-1,j}$ belongs  to $t\, \qs(T)$, and
 is thus algebraic.
More generally, for $i\ge 1$ and $j\ge 0$
 having the same parity, 
the series $C_{-i,j}$ is D-finite, of the form
 $$
- \frac 1 3 Q_{i-2,j} + t^{\min(i,j)} \Rat(Z)
$$
for some rational function $\Rat$. It is  transcendental
as soon as $i\ge 2$.

Finally, for $i\ge 0$ and $j\ge 0$ having the same parity,  the series $C_{i,j}$  is of the form
 $$
 \frac 1 3 Q_{i,j} + t^{\min(i,j)} \Rat(Z).
$$
It is D-finite and transcendental.
\end{cor}
Here are some  examples:
\beq\label{coincidence2}
tC_{-1,1}=\frac{(T-1)(11+6T-T^2)}{(T+3)^3}
,
\eeq
$$
C_{-2,0}= -\frac 1  3 Q_{0,0}+ \frac {32}3\, \frac{Z^3(1+Z+3Z^2-Z^3)}{(Z+1)(Z^2+3)^3},
$$
$$
C_ {0,0}= \frac 1 3 Q_{0,0} + \frac {64}3\, \frac{Z^3(1+Z+3Z^2-Z^3)}{(Z+1)(Z^2+3)^3}.
$$
The similarity between the last two expressions comes
from~\eqref{PM-diag}. Finally, we have found closed form expressions
for walks ending on the boundary of the cone $\mathcal C$. For instance:
\begin{align}
  c_{0,0}(2n)&=\frac{16^n}9 \left(  3 \, \frac{(1/2)_n^2}{(2)_n^2}
+  8 \,  \frac{(1/2)_n(7/6)_n}{(2)_n(4/3)_n}
-  2 \,  \frac{(1/2)_n(5/6)_n}{(2)_n(5/3)_n}\right), \nonumber
\\
c_{-2,0}(2n)&=\frac{16^n}9 \left(-  3 \, \frac{(1/2)_n^2}{(2)_n^2}
+  4 \,  \frac{(1/2)_n(7/6)_n}{(2)_n(4/3)_n}
-    \frac{(1/2)_n(5/6)_n}{(2)_n(5/3)_n}\right), \label{expressions-diag}
\\
  c_{-4,0}(2n)&=\frac{16^n}{3^5} \left(-3^5 n \, \frac{(1/2)_n^2}{(2)_n(2)_{n+1}}
+  4 ({21n^2+30n-14})\,
\frac{(1/2)_n(7/6)_n}{(2)_{n+1}(4/3)_{n+1}}\right.  \nonumber
\\&\left.
\hskip 50mm -    7 ({3n^2+3n-10}) \,
\frac{(1/2)_n(5/6)_n}{(2)_{n+1}(5/3)_{n+1}}\right).  \nonumber
\end{align}
It seems that this pattern persists, that is, that $c_{-2i,0}$ is a sum of three hypergeometric terms (we have checked this for $0\le i \le 4$.)
There is no such expression for walks ending at  $(-3,1)$, $(-1,1)$, $(1,1)$, $(-2,2)$, nor $(0,2)$.

\bigskip
Our starting point is of course the functional equation obtained by
constructing walks recursively. It reads
$$
K(x,y)\C(x,y)=1 -t\by(x+\bx) \Cn(\bx) -t\bx(y+\by) \Cn(\by)  -t\bx\by \C_{0,0},
$$
where $\Cn(\bx)$ is still given  by~\eqref{Tm-def}, the kernel is 
$
K(x,y)=1-t(x+\bx)(y+\by),
$
and  $C_{0,0}$ is the coefficient of $x^0y^0$ in $C(x,y)$. As in the
square lattice case, the kernel is invariant by the transformations
$x\mapsto \bx$ and $y\mapsto \by$.

\subsection{Reduction to an equation with orbit sum zero}
\label{sec:diag-os0}

Let us compare this  equation to the one that describes  walks confined to $\mathcal Q$:
$$
K(x,y) Q(x,y)=1-t\by(x+\bx) \Qp(x) - t\bx(y+\by) \Qp(y)  +t\bx\by Q_{0,0},
$$
with $Q_{0,0}:=Q(0,0)$ and $Q_+(x)=Q(x,0)$. The  orbit
equations of $Q$ and $C$ are still given by~\eqref{orbit-Q} and~\eqref{orbit-C}, and in particular they have the same right-hand side $(x-\bx)(y-\by)$.
This leads us to introduce a series $A(x,y)$ defined
by~\eqref{A-def}, again with $a=1/3$. The
equation satisfied by $A$ (the counterpart of~\eqref{eq-A}) is 
$$
K(x,y)A(x,y)= (2+\bx^2+\by^2)/3-t\by(x+\bx)A_-(\bx)
-t\bx(y+\by)A_-(\by)-t\bx\by A_{0,0}.
$$
The corresponding orbit sum is of course zero.

\subsection{Reduction to a quadrant-like problem}

In the series $A$, we separate the contributions of the three quadrants
by introducing the series $P$ and $M$ given by~\eqref{A-expr}. Given
that the orbit sum of $A$ is zero, these two series are 
still related by~\eqref{PM-bis}. We now follow the lines of Section~\ref{sec:quadrant-red}
to obtain the quadrant-like equation~\eqref{func-M-diag} for
$M(x,y)$. It is  the diagonal counterpart of~\eqref{eqcat-M}.

\subsection{Cancelling the kernel: an equation between
  $\boldsymbol{M(0,x)}$,  $\boldsymbol{M(0, \bx)}$ and
  $\boldsymbol{M(x,0)}$}
\label{sec:cancel-diag}
As a polynomial in $y$, the kernel $K$ admits only one root in the
ring of formal series in $t$:
\beq\label{Y-sol-diag}
Y= \frac{1-\sqrt{1-4t^2(x+\bx)^2}}{2t(x+\bx)} = (x+\bx) t +
(x+\bx)^3t^3+O(t^5). 
\eeq
We follow the steps  of Section~\ref{sec:cancel} to obtain the counterpart
of~\eqref{catM}:
\beq\label{catM-diag}
(x+\bx)(Y-1/Y) \left( x M(0,x)-2 \bx M(0,\bx)\right) -2\bx Y/t+3(x+\bx) M(x,0)+3M_{1,0}=0.
\eeq

\subsection{An equation between $\boldsymbol{M(0,x)}$ and $\boldsymbol{M(0,\bx)}$}
 The discriminant occurring in the expression~\eqref{Y-sol-diag} of $Y$ is 
$$1-4t^2(x+\bx)^2=1-4t^2(x^2+1)(\bx^2+1).$$
 It is an even function of $x$. The same
holds for the series $xM(x,0)$ and $xM(0,x)$. This suggests to define:
$$
\Delta= 1-4t^2(x+1)(\bx+1),
$$
\beq\label{RS-def-diag}
R(x)= t^2/\sqrt x\, M(\sqrt x,0), \qquad 
S(x)=t\sqrt x\,  M(0,\sqrt x).
\eeq
Then $R(x)$ and $S(x)$ both belong to $\qs[x][[t^2]]$, and~\eqref{catM-diag} gives:
\beq\label{eqRS-diag}
\sqrt{\Delta(x)}\left( S(x)-2S(\bx) -\frac 1 {1+x} \right) = 3(x+1)R(x)+3R_0-\frac 1
{1+x},
\eeq
with $R_0=R(0)=t^2M_{1,0}$. Expanding around $x=-1$ gives
\beq\label{R0-Sm1}
3R_0+S(-1)=0,
\eeq
which will be useful later.

As in Section~\ref{sec:M0x-M0bx}, we square~\eqref{eqRS-diag} and then extract the negative part. This gives
$$
\Delta(x) \left( S(\bx)^2+ \frac{S(\bx)}{1+x}\right) -[x^<]\left(
  \Delta(x) S(x) S(\bx)\right)= \bx t^2 + \frac{S(-1)}{x+1}.
$$
As before, we denote by $P_0$ the coefficient of $x^0$ in $\Delta(x)
S(x) S(\bx)$. Reconstructing this series finally gives the counterpart
of~\eqref{eqcat-S}:
$$
\Delta(x) \left( S(x)^2+S(\bx)^2-S(x)S(\bx) + \frac{S(x)}{\bx +1} +
  \frac{S(\bx)}{x+1}\right) = S(-1)-P_0+ t^2(x+\bx).
$$

\subsection{An equation for  $\boldsymbol{M(0,x)}$ only}
\label{sec:M0x-diag}
We multiply the previous equation by $S(x)+S(\bx)+1$. The non-negative
part of the resulting equation reads:
\begin{multline*}
   \Delta(x)\left( S(x)^3+ \frac{2x+1}{x+1}S(x)^2+\frac{S(x)}{\bx+1}\right)
=\\
 \left( S(-1)-P_0+  t^2(x+\bx)\right) (S(x)+1)+t^2S_1-
\frac{S(-1)(S(-1)+1)}{x+1}-t^2\bx
\end{multline*}
with $S_1:=S'(0)$. Extracting the constant term in $x$ gives
$P_0+S(-1)^2=2t^2S_1$, which allows us to rewrite the above equation as
\begin{multline}\label{eqcatS-diag}
   \Delta(x)\left( S(x)^3+ \frac{2x+1}{x+1}S(x)^2+\frac{S(x)}{\bx+1}\right)
=\\
 \left(   t^2(x+\bx)-F_0\right) (S(x)+1)+t^2S_1-
\frac{2t^2S_1-F_0}{x+1}-t^2\bx,
\end{multline}
with $F_0=P_0-S(-1)$.

\subsection{The generalized quadratic method: algebraicity of $\boldsymbol{A(x,y)}$}
\label{sec:quad-diag}
We now apply the generalized quadratic method of
Section~\ref{sec:quad}. We denote
\begin{multline*}
  \Pol(x_0,x_1,x_2,t,x)=   \Delta(x)\left( x_0^3+ \frac{2x+1}{x+1}x_0^2
+\frac{x_0}{\bx+1}\right)
\\
- \left(   t^2(x+\bx)-x_2\right) (x_0+1)-t^2x_1+\frac{2t^2x_1-x_2}{x+1}+t^2\bx.
\end{multline*}
Then~\eqref{quad1} holds with $A_1=S_1$ and $A_2=F_0$. We find
that~\eqref{quad2} now admits two solutions $X_0$ and $X_1$. Computing
their first few coefficients leads us to  conjecture that
$$
X_0= \frac{1-2t-\sqrt{1-4t}}{2t}, \qquad X_1=
-\frac{1+2t-\sqrt{1+4t}}{2t}.
$$
This is proved by eliminating $x_1$ and $x_2$ between $\Pol$,
$\Pol_{x_0}$ and $\Pol_x$, since
\begin{multline*}
  2 \Pol-(1+2x_0)\Pol_{x_0}-(x^2-1) \Pol_x=\\
-
\left(1-2t - t(x+\bx)\right)\left(1+2t + t(x+\bx)\right)\frac{(1+2x_0)^2(x+x_0(1+x))}{x+1}.
\end{multline*}
The series $X_0$ and $X_1$ cancel the first and second factor, respectively.

It remains to say that the discriminant of $\Pol(x_0, A_1,A_2,t,x)$
with respect to $x_0$  admits roots at $x=X_0$ and $x=X_1$. This
gives two polynomial relations between $A_1$ and $A_2$, that is,
between $S_1$ and $F_0$, from which we derive:
\begin{multline*}
  19683{t}^{6}S_1 ^{4}+2187{t}^{4} \left( 20{t}^{2}-1
 \right) S_1 ^{3}+81{t}^{2} \left( 11{t}^{2}-1 \right) 
 \left( 38{t}^{2}-1 \right) S_1  ^{2}\\
+ \left( 92{t}^{2}-1
 \right)  \left( 11 t^{2}-1 \right) ^{2}S_1 +{t}^{2} \left( 
1331{t}^{4}-107{t}^{2}+1 \right) =0
\end{multline*}
and
\begin{multline*}
  27 F_0^{4}+27  \left( 8 {t}^{2}-1 \right) F_0^{3}+9  \left( 2
    t+1 \right)  \left( 2 t-1 \right)  \left( 10 {t}^{2}-1 \right)
  F_0^{2}\\
+ \left( 224 {t}^{6}-68 {t}^{4}+16 {t}^{2}-1 \right) F_0+{t}^{2} \left( 48 {t}^{6}+88 {t}^{4}-20 {t}^{2}+1 \right)=0.
\end{multline*}
Note that the equation satisfied by $S_1$  is the same  as in 
the square lattice case (see~\eqref{S1-eq-sq}).

From this point on, we conclude that the series $S(x)$ (or $M(0,x)$),
$R(x)$ (or $M(x,0)$), $M(x,y)$ and finally $P(x,y)$ and $A(x,y)$ are
algebraic, using, in this order~\eqref{eqcatS-diag},
\eqref{eqRS-diag}, \eqref{func-M-diag}, \eqref{PM-diag} and~\eqref{A-expr}.

 \subsection{Rational parametrization and degrees}
The end of the proof of Theorem~\ref{thm:diag} is very similar to
Section~\ref{sec:param}. The above equations for $S_1$ and $F_0$ factor when
$t^2$ is parametrized by $T$, and one obtains
$$
S_1= \frac{(T-1)(11+6T-T^2)}{(3+T)^3},
$$
while
$$
F_0= \frac{(1-T)(3T^3-29T^2-15T+9)}{128 T^3}.
$$
We plug these expressions in the equation~\eqref{eqcatS-diag}
defining $S(x)$, parametrize $x$ by the series $V$ defined by~\eqref{V-def}, and obtain a rational expression of $S(x)$ in terms of
$T$ and~$V$, which is equivalent to the expression~\eqref{M0x-diag} of
$t\sqrt x\, M(0,\sqrt x)$.

We then consider~\eqref{eqRS-diag}. Our first task is to determine
$R_0$, or equivalently $S(-1)$ (see~\eqref{R0-Sm1}). In order to do so, we
replace $x$ by $-1$ in the cubic equation defining $S(x)$ over
$\qs(x,T)$. The resulting equation factors once we write $T=Z^2$,
giving
$$
R_0=-\frac 1 3 \, S(-1)={\frac { \left( Z-1 \right)  \left(-
      {Z}^{3}+3\,{Z}^{2}+Z+1 
 \right) }{24{Z}^{3}}}.
$$
Then, the term
$D(x)=S(x)-2S(\bx)-1/(1+x)$ is found to be bicubic over
$\qs(T,x)$. Combined with the above value of $R_0$, this gives for
$R(x)$ an equation of degree 6
over $\qs(Z,x)$, which factors into two cubic terms, and factors even
further when parametrizing $x$ by $T$ (or  $Z$) and $V$. This gives a rational
expression for $R(x)$ in terms of $Z$ and $V$, which is equivalent to
the expression~\eqref{Mx0-diag} of $t/\sqrt x\, M(\sqrt x,0)$.

It remains to prove that $M(x,y)$ and $P(x,y)$ have degree 72 over
$\qs(x,y,t)$. We have proved that $M(x,0)$ and $M(0,x)$ belong to
$\qs(t, x,Z, V_2)$, where $V_2$ denotes the series $V$ with $x$ replaced
by $x^2$. The functional equations~\eqref{func-M-diag} and~\eqref{PM-diag} defining
$M(x,y)$ and $P(x,y)$ prove that both series belong to $\qs(t,x,y,Z, V_2,
\tilde V_2)$, where $\tilde V_2$ is $V_2$ with $x$ replaced by
$y$. Since $V_2$ is cubic over $\qs(x,Z)$, it follows that $M(x,y)$
and $P(x,y)$ have degree at most 72 over $\qs(x,y,t)$. Computing, by
successive eliminations, their minimal equation when $x=2$ and $y=3$
shows that this bound is tight.

\subsection{Walks ending at a prescribed position}
The argument is similar to that of Section~\ref{sec:prescribed}. The result boils
down to proving that $R(x)$, seen as a series in $x$, has coefficients
in $\qs(Z)$, while $S(x)$
has coefficients in $\qs(T)$.  

The (cubic) equation over $\qs(x,Z)$ satisfied by $R(x)$ has degree 1
in $R$ when $x=0$, and this permits a recursive computation of the
coefficients of $R(x)$ in $\qs(Z)$. A similar statement holds for the
equation over $\qs(x,T)$ satisfied by $S(x)/x$.

\section{Starting at $\boldsymbol{(-1,0)}$ on the square lattice}
\label{sec:square-asym}
We now return to the ordinary square lattice, but change the starting
point to $(-1,0)$. The
$x/y$ symmetry is lost, which complicates the derivation a
bit.  On the other hand, the series $\C(x,y)$ counting walks confined
to $\mathcal C$ now satisfies a functional 
equation with orbit sum zero, and turns out to be algebraic.
That algebraicity is sensitive to the starting point is not a new
phenomenon: for instance, quadrant walks with steps
$\nearrow,\downarrow,\leftarrow$, known as Kreweras' walks, are algebraic when
starting at $(0,0)$~\cite{gessel-proba}, but transcendental when starting at
$(1,0)$~\cite{kilian-private}. 

\begin{theorem}\label{thm:square-asym}
  The \gf\ of square lattice walks starting at $(-1,0)$ and confined
  to the cone $\mathcal C$ is algebraic.  Let $P(x,y)$ (resp. $\bx
  L(\bx,y)$, $\by B(x,\by)$) denote the \gf\ of such walks ending in the
  first quadrant (resp. at a negative abscissa, at a negative
  ordinate). These three series are algebraic of degree $72$ over $\qs(x,y,t)$.

More precisely, $P$ can be expressed in terms of $L$ and $B$ by:
\beq\label{P-LB}
P(x,y)= \bx\left( L(x,y)-L(0,y)\right) + \by
\left(B(x,y)-B(x,0)\right),
\eeq
the series $L$ and $B$ satisfy
\begin{multline}
  \label{eqL}
\left( 1-t(x+\bx +y+\by)\right)  \left(2L(x,y)-L(  0,y)\right)=\\
1-2t\by L(x,0)+t(x -\bx)
L(0,y)+t\by B(0,y)+t\by L_{0,0} -t\by B_{0,0}
\end{multline}
and
\begin{multline}
  \label{eqB}
\left( 1-t(x+\bx +y+\by)\right)  \left(2B(x,y)-B(x,0)\right)=\\
t(y-\by) B(x,0)-2t\bx B(0,y)+t\bx
L(x,0)
-t\bx L_{0,0} +t\bx B_{0,0},
\end{multline}
where $L_{0,0}=L(0,0)$ and $B_{0,0}=B(0,0)$, and each specialization
$L(x,0)$, $L(0,x)$, $B(x,0)$ and $B(0,x)$ has  degree $24$ over $\qs(t,x)$.

Moreover, these algebraic series have rational
parametrizations. Defining the series $T$, $Z$ and $U$ as in
Theorem~{\rm\ref{thm:square}}, the series $L(xt,0)$, $L(0,xt)$, $B(xt,0)$,
$B(0,xt)$ admit rational expressions in terms of $Z$ and $U$, given in
 Appendix~\ref{sec:param-sq-asym}.
\end{theorem}
\begin{cor}[\bf Walks ending at a prescribed position]\label{cor:coeffs-square-asym}
  Let $T$ be the unique series in $t$ with constant term $1$ 
  satisfying~\eqref{T-def}, and let  $Z=\sqrt T $.

 For any $(i,j)$ in $\mathcal C$,  the series $C_{i,j}$ belongs  to $t^{i+j-1} \qs(Z)$, and
 is thus algebraic. 
\end{cor}
\noindent Sometimes $C_{i,j}$ even belongs to  $t^{i+j-1} \qs(T)$. Here are some examples:
$$
tC_{0,0}=-{\frac { \left( T-1 \right)  \left( T^{2}-6\,T-11 \right) }
{ \left( T+3 \right) ^{3}}}, \qquad tC_{-2,0}=16\,{\frac { \left( T-1 \right)  }{ \left( T+3 \right) ^{3}}},
$$
$$
tC_{0,-2}={\frac { \left( T-1 \right)^2(5-T) }{ \left( T+3 \right) ^{3}}}
, \qquad C_{-1,0}=64\,{\frac {Z^{3}}{ \left( Z^{2}+3 \right) ^{3}}},
$$
$$
tC_{-1,1}= -16\,{\frac {Z^{2} \left( Z-1 \right)  \left( Z-3 \right)
  }{ \left( Z^{2}+3 \right) ^{3}}}, \qquad
C_{0,-1}=-32\,{\frac {Z^{3} \left( Z-1 \right)  \left( Z^{2}-2\,Z-1 \right) }{ \left( Z+1 \right)  \left( Z^{2}+3 \right) ^{3}}}.
$$
From these expressions, we can look for hypergeometric forms of the
coefficients.  In this way, we find
$$
c_{-1,0}(2n)= \frac{16^n}{3}\left( 2\, \frac {(1/2)_n (7/6)_n}{(2)_n(4/3)_n}+\frac{(1/2)_n (5/6)_n}{(2)_n (5/3)_n}\right),
$$
$$
c_{0,-1}(2n) =\frac{2 \cdot 16^n}{3}\left(  \frac {(1/2)_n (7/6)_n}{(2)_n(4/3)_n}
-\frac{(1/2)_n (5/6)_n}{(2)_n (5/3)_n}\right),
$$
but no other simple expression in the vicinity of the origin. The
reflection principle directly relates the number of Gessel walks
ending at $(0,0)$ to the above numbers (see Section~\ref{sec:gessel}).

\subsection{Reduction to two quadrant-like problems}
Let $\C(x,y)$ be  the \gf\ of square lattice walks starting at $(-1,0)$ and confined
to  $\mathcal C$. It satisfies the functional equation
\beq\label{eqC-sq-asym}
K(x,y)C(x,y)=\bx  -t\by \Ch(\bx) -t\bx \Cv(\by),
\eeq
where $K(x,y)=1-t(x+\bx+y+\by)$ is the kernel,
\beq\label{Ch}
\Ch(\bx)= \sum_{i<0,n\ge 0}c_{i,0}(n) x^i t^n
\eeq
and
\beq\label{Cv}
\Cv(\by)= \sum_{j<0,n\ge 0}c_{0,j}(n) y^j t^n.
\eeq
Due to the constant term $\bx$ in~\eqref{eqC-sq-asym} (instead of $1$ in Section~\ref{sec:square}), the orbit sum vanishes:
\beq\label{oe-square-asym}
xyC(x,y)-\bx y C(\bx ,y) +\bx \by \C(\bx,\by)-x\by C(x,\by)=0.
\eeq
We write
$$
C(x,y)= P(x,y)+\bx L(\bx,y)+ \by B(x,\by),
$$
where $P(x,y), L(x,y)$ and $B(x,y)$ belong to $\qs[x,y][[t]]$ (the letter $P$ stands for \emm positive,, and the letters $L$ and $B$
for \emm left, and \emm below,, respectively).
We plug this expression of $C$ in the orbit
equation~\eqref{oe-square-asym}, and extract the positive part in
$x$ and~$y$. This gives the expression~\eqref{P-LB} of $P$ in terms of $L$ and
$B$. We can thus express $C$ in terms of $L$ and $B$ as well.

We now return to the equation~\eqref{eqC-sq-asym} defining $C$, and 
replace $C$ by its expression in terms of $L$ and $B$. Extracting from
the resulting equation the negative part in $x$ gives
$$
\bx K(x,y) L(\bx,y)=-t\bx\by L(\bx,0)-tL(0,y)+t\bx \by B(0,y) +t\bx
L_x(0,y)-t\bx\by B_{0,0}+\bx,
$$
with $L_x=\partial L/\partial x$. Extracting from this the coefficient
of $\bx$ gives an expression of $L_x(0,y)$ in terms of $L(0,y)$ and
$B(0,y)$, which, plugged in the previous equation, leads to:
\begin{multline*}
  \bx K(x ,y)\left(L(\bx ,y)-L(0,y)/2\right)=\bx /2-t\bx \by
L(\bx ,0)\\
+t(\bx ^2-1)L(0,y)/2+t\bx \by B(0,y)/2 +t\bx \by L_{0,0}/2-t\bx \by B_{0,0}/2,
\end{multline*}
which is equivalent to~\eqref{eqL}, upon replacing $x$ by $\bx$.

Repeating the procedure with  $y$ instead of $x$ leads to
the equation~\eqref{eqB} satisfied by $B$.

Equations~\eqref{eqL} and~\eqref{eqB} both involve the series $L$ and
$B$. To obtain two decoupled equations, we define
\beq\label{MN-def}
M(x,y)=L(x,y)+B(y,x) \qquad \hbox{and} \qquad N(x,y)=L(x,y)-B(y,x).
\eeq
Then
$$
K(x,y)\left(2M(x,y)-M(0,y)\right)=1-2t\by M(x,0)+t(x-\bx) M(0,y)+t\by
M(y,0)
$$
and
$$
K(x,y)\left(2N(x,y)-N(0,y)\right)=1-2t\by N(x,0)+t(x-\bx)N(0,y)-t\by
N(y,0)+2t\by N_{0,0}.
$$
These equations are extremely close to~\eqref{eqcat-M}, and we solve them in
exactly the same way. Below we give a few details on some steps of the
procedure, but we otherwise refer to the Maple session available
on the author's \href{http://www.labri.fr/perso/bousquet/publis.html}{webpage}.
Having found $M(x,y)$ and $N(x,y)$,
we  reconstruct $L$ and $B$ thanks to~\eqref{MN-def}. The
remaining results (degrees, nature of the coefficients)  are
established as in Section~\ref{sec:square}.

\subsection{Solving the equation for $\boldsymbol{M(x,y)}$}
We introduce the series $R$ and $S$ related to $M$
by~\eqref{RS-def}. Following the steps of Sections~\ref{sec:cancel}
to~\ref{sec:M0x} leads to a cubic
equation for $S(x)$ 
 which, as~\eqref{eqcatS}, involves two
additional unknown series, namely $R_0=S_1$ and $S_2$.

We apply to this equation the generalized quadratic method of
Section~\ref{sec:quad}. Two series cancel $\Pol_{x_0}$. One of them is
$$
X_0= \frac{1-\sqrt{1-16t^2}}{4t},
$$
and the other satisfies $X_1=X$ with
$$
(2t+X+1/X-t(X^2+1/X^2))S(X)(2+S(X))+t+(X+1/X)-t(X^2+1/X^2)/2=0.
$$
Proceeding exactly as in Section~\ref{sec:quad}, we derive from this
that  $S_1$ and $S_2$ have degree~8 over $\qs(t)$. After introducing
the series $Z$, we obtain for  $S_1/t$ and
 $S_2$  rational expressions in terms of  $Z$.

Then one finds that $S(xt)=xt^2M(0,xt)$ (which is an even function of $t$) has degree 24
over $\qs(x,t)$, and is cubic over $\qs(Z,x)$. It can be expressed
rationally in terms of the series $Z$ and $U$.

For the series $R(xt)/t$, the degree is only 12 over $\qs(x,t)$, and 3
over $\qs(T,x)$. Again, introducing $U$ factors the equation and gives
a rational expression of $R(xt)/t=M(xt,0)$ in terms of $T$ and $U$.

\subsection{Solving the equation for $\boldsymbol{N(x,y)}$}
We introduce  series $R$ and $S$ related to $N$ in the same way they
were related to $M$ before
(see~\eqref{RS-def}). Following the steps of Sections~\ref{sec:cancel}
to~\ref{sec:M0x} leads to a cubic
equation for $S(x)$,  which now involves only one
additional unknown series, namely $R_0=S_1$.

We apply to this equation the generalized quadratic method of
Section~\ref{sec:quad}. One  series $X$ cancels $\Pol_{x_0}$, and we derive from
its existence an equation of degree 8 for $S_1$ over $\qs(t)$. Again, $S_1/t$ has a rational expression in terms of the series $Z$.

Then one finds that $S(xt)=xt^2N(0,xt)$ (which is an even function of $t$) has degree 12
over $\qs(x,t)$, and is cubic over $\qs(T,x)$. It can be expressed
rationally in terms of the series $T$ and $U$.

For the series $R(xt)/t$, the degree is 24 over $\qs(x,t)$, and 3
over $\qs(Z,x)$.  Introducing $U$ factors the equation and gives
a rational expression of $R(xt)/t=N(xt,0)$ in terms of $Z$ and $U$.

\section{Starting at $\boldsymbol{(-2,0)}$ on the diagonal square lattice}
\label{sec:diag-asym}
We finally return to the diagonal lattice and change the starting
point to $(-2,0)$. The orbit sum associated with the \gf\ $C(x,y)$ is
non-zero, and the \gf\ $Q(x,y)$ 
counting quadrant walks with diagonal steps, given
by~\eqref{Q-diag},  enters the picture
again. The \gf\ for walks starting from $(-2,0)$ and confined to the
cone $\cC$ now differs from
$$
-\frac 1  3  \, \left( Q(x,y)-\bx^2Q(\bx,y)-\by ^2Q(x,\by)\right)
$$
by an algebraic series.
Comparing with Theorem~\ref{thm:diag} shows that the above D-finite part is
the \emm opposite, of what it was when starting from $(0,0)$.

\begin{theorem}
  \label{thm:diag-asym}
   The \gf\ of walks on the diagonal
  square lattice that start from $(-2,0)$, remain in $\cC$, and end in the
  first quadrant (resp. at a negative abscissa, at a negative
  ordinate) is
$$
-\frac 1 3 \, Q(x,y) +P(x,y)
\qquad \left(\hbox{resp. } \frac 1 3 \, \bx^2 Q(\bx,y) +\bx L(\bx,y), \qquad
\frac  1 3 \, \by^2 Q(x,\by) +\by B(x,\by)\right),
$$
where $P(x,y)$, $L(x,y)$ and $B(x,y)$  are algebraic of degree $72$
over $\qs(x,y,t)$. 

More precisely, $P$ can be expressed in terms of $L$ and $B$ by:
\beq\label{P-LB-diag}
P(x,y)= \bx\left( L(x,y)-L(0,y)\right) + \by
\left(B(x,y)-B(x,0)\right),
\eeq
the series $L$ and $B$ satisfy
\begin{multline}
  \label{eqL-diag}
\left( 1-t(x+\bx)(y+\by)\right)  \left(2L(x,y)-L(  0,y)\right)=\\
4x/3-2t\by (x+\bx)L(x,0)+t(x -\bx)(y+\by)
L(0,y)+t(1+\by^2) B(0,y) -t\by B_{0,1}
\end{multline}
and
\begin{multline}
  \label{eqB-diag}
\left( 1-t(x+\bx)(y+\by)\right)  \left(2B(x,y)-B(x,0)\right)=\\
-2y/3 +t(x+\bx)(y-\by) B(x,0)-2t\bx(y+\by) B(0,y)+t(1+\bx^2)
L(x,0)
-t\bx L_{1,0},
\end{multline}
where $L_{1,0}$ (resp. $B_{0,1}$) denotes the coefficient of $x^1y^0$
(resp. $x ^0 y^1$) in $L(x,y)$
(resp. $B(x,y)$).  The specializations $L(0,x)$ and  $B(x,0)$   have
 degree $12$ over $\qs(t,x)$, while  $L(x,0)$ and $B(0,x)$ have degree $24$.

Moreover, these algebraic series have rational
parametrizations. Defining the series $T$, $Z$ and $V$ as in
Theorem~{\rm\ref{thm:diag}}, the series $\sqrt x L(\sqrt x,0)$,
$t\sqrt x L(0,\sqrt x)$, $t\sqrt x B(\sqrt x,0)$, and
$\sqrt x B(0,\sqrt x)$ (which belong to $\qs[x][[t^2]]$) admit
rational expressions in terms of $Z$ and $V$, given in 
 Appendix~\ref{sec:param-diag-asym}.
\end{theorem}

\begin{cor}[\bf Walks ending at a prescribed position]\label{cor:coeffs-diag-asym}
 Let $T$ be the unique series in $t$ with constant term $1$
  satisfying~\eqref{T-def}, and let  $Z=\sqrt T $.

 For $j\ge 0$, the series $C_{-1,j}$ and $C_{j,-1}$ belong  to $t\, \qs(T)$, and
 are thus algebraic.
More generally, for $i\ge 1$ and $j\ge 0$
 having the same parity, 
the series $C_{-i,j}$ and  $C_{j, -i}$ are D-finite, of the form
 $$
 \frac 1 3\, Q_{i-2,j} + t^{\min(i,j)} \Rat(Z) 
$$
for some rational function $\Rat$. They are  transcendental
as soon as $i\ge 2$.

Finally, for $i\ge 0$ and $j\ge 0$, having the same parity,  the series $C_{i,j}$  is of the form
 $$
 -\frac 1 3\, Q_{i,j} + t^{\min(i,j)} \Rat(Z).
$$
It is D-finite and transcendental.
\end{cor}
Here are some examples:
$$
tC_{-1,1}= \frac{16(T-1)}{(T+3)^3},\qquad
tC_{-1,3}= \frac{64(T-1)^2(T+1)(7-T)}{(T+3)^6},
$$
$$
C_{-2,0}=\frac 1 3 Q_{0,0}+ \frac {32Z^3(5+5Z-3Z^2+Z^3)}{3(Z+1)(Z^2+3)^3},
$$
$$
C_{0,0}=-\frac 1 3 Q_{0,0}+ \frac {32Z^3(1+Z+3Z^2-Z^3)}{3(Z+1)(Z^2+3)^3},
$$
$$
C_{0,-2}= \frac 1 3 Q_{0,0}- \frac {64Z^3(2+2Z-3Z^2+Z^3)}{3(Z+1)(Z^2+3)^3},
\qquad
tC_{1,-1}= \frac{(T-1)^2(5-T)}{(T+3)^3}.
$$
We found hypergeometric expressions for the number of $n$-step walks
starting from $(-2,0)$ and ending at a prescribed point of the
boundary of $\cC$:
\begin{align*}
  c_{0,0}(2n)&=\frac{16^n}9 \left(-  3 \, \frac{(1/2)_n^2}{(2)_n^2}
+  4 \,  \frac{(1/2)_n(7/6)_n}{(2)_n(4/3)_n}
-    \frac{(1/2)_n(5/6)_n}{(2)_n(5/3)_n}\right),
\\
c_{-2,0}(2n)&=\frac{16^n}9 \left(  3 \, \frac{(1/2)_n^2}{(2)_n^2}
+  2 \,  \frac{(1/2)_n(7/6)_n}{(2)_n(4/3)_n}
+4    \frac{(1/2)_n(5/6)_n}{(2)_n(5/3)_n}\right),
\\
c_{0,-2}(2n)&=\frac{16^n}9 \left(  2 \, \frac{(1/2)_n^2}{(2)_n^2}
-5 \,  \frac{(1/2)_n(7/6)_n}{(2)_n(4/3)_n}
+2    \frac{(1/2)_n(5/6)_n}{(2)_n(5/3)_n}\right),
\\
  c_{-4,0}(2n)&=\frac{16^n}{3^5} \left(  3^5 n\, \frac{(1/2)_n^2}{(2)_n(2)_{n+1}}
+  2 (21n^2+30n+16)\,
\frac{(1/2)_n(7/6)_n}{(2)_{n+1}(4/3)_{n+1}}\right.
\\ & \hskip 40mm
\left. +4 (39n^2+66n-10)   \frac{(1/2)_n(5/6)_n}{(2)_{n+1}(5/3)_{n+1}}\right).
\end{align*}
This pattern persists at least up to $c_{-8,0}(2n)$ and
$c_{0,-8}(2n)$.

The proof of these results  combines those of the last two
sections: Section~\ref{sec:diag}, which dealt with walks starting from
$(0,0)$ on the diagonal lattice, and Section~\ref{sec:square-asym},
which dealt with walks starting at $(-1,0)$ on the square
lattice. There is one new difficulty in the application of the generalized
quadratic method, because none of the auxiliary series $X_i$ is easy
to guess. We explain how to overcome this problem, but otherwise
simply write down some intermediate results of the derivation.

\subsection{Reduction to two quadrant-like problems}
Let $\C(x,y)$ be  the \gf\ of walks starting at $(-2,0)$ and confined
to  the cone $\mathcal C$ in the diagonal square lattice. It satisfies the functional equation
$$
K(x,y)C(x,y)=\bx^2  -t\by(x+\bx) \Ch(\bx) -t\bx(y+\by)
\Cv(\by)-t\bx\by C_{0,0},
$$
where $K(x,y)=1-t(x+\bx)(y+\by)$ is the kernel, and the series
$\Ch(\bx)$ and $\Cv(\by)$ are defined by~(\ref{Ch}-\ref{Cv}).
The orbit sum is now $-(x-\bx)(y-\by)$. This is the opposite of the
orbit sum for quadrant walks (see Section~\ref{sec:diag-os0}), and this leads us
to introduce the series 
$$
A(x,y):=C(x,y)+ \frac 1  3 \left( Q(x,y)-\bx^2Q(\bx,y)-\by ^2Q(x,\by)\right).
$$
The equation satisfied by $A$ reads
\beq\label{eqA-diag-asym}
K(x,y)A(x,y)=\frac 1 3 \left(1+2\bx^2-\by^2 \right) -t\by(x+\bx) \Ah(\bx)
-t\bx(y+\by)\Av(\by) -t\bx\by A_{0,0},
\eeq
and now the orbit sum vanishes. 
We write
$$
A(x,y)= P(x,y)+\bx L(\bx,y)+ \by B(x,\by),
$$
where $P(x,y), L(x,y)$ and $B(x,y)$ belong to $\qs[x,y][[t]]$.
We plug this expression of $A$ in the above
equation, and extract the positive part in
$x$ and $y$. This gives the expression~\eqref{P-LB-diag} of $P$ in terms of $L$ and
$B$. We can thus express $A$ in terms of $L$ and $B$ as well.

We plug this expression of $A$ in the equation~\eqref{eqA-diag-asym}.  Extracting  the negative part in $x$ gives an equation
which is equivalent to~\eqref{eqL-diag}, upon replacing $x$ by $\bx$.

Symmetrically, extracting  the negative part in $y$ leads to
the equation~\eqref{eqB-diag} satisfied by $B$.

 As in the square
lattice case, we can  decouple the series $L$ and $B$ by defining
\beq\label{MN-def-diag}
M(x,y)=L(x,y)+B(y,x) \qquad \hbox{and} \qquad N(x,y)=L(x,y)-B(y,x).
\eeq
Then
\begin{multline*}
  K(x,y)\left(2M(x,y)-M(0,y)\right)=2x/3-2t\by(x+\bx)
M(x,0)\\
+t(x-\bx)(y+\by) M(0,y)+t(1+\by^2)M(y,0)- t\by M_{1,0}
\end{multline*}
and
\begin{multline*}
  K(x,y)\left(2N(x,y)-N(0,y)\right)=2x-2t\by(x+\bx) N(x,0)\\
+t(x-\bx)(y+\by)N(0,y)-t(1+\by^2)N(y,0)+t\by N_{1,0}.
\end{multline*}
The first  equation is exactly the one we met when counting walks
starting at $(0,0)$ on the diagonal lattice
(see~\eqref{func-M-diag}), and we  only have to solve the other one. Below we give a few details on some steps of the
procedure, but we again refer to the Maple session available
on the author's \href{http://www.labri.fr/perso/bousquet/publis.html}{webpage}.
Having found $M(x,y)$ and $N(x,y)$,
we  reconstruct $L$ and $B$ thanks to~\eqref{MN-def-diag}. The
remaining results (degrees, nature of the coefficients)  are
established as in Section~\ref{sec:diag}.

\subsection{Solving the equation for $\boldsymbol{N(x,y)}$}
We introduce  series $R$ and $S$ related to $N$ in the same way they
were related to $M$ in~\eqref{RS-def-diag}. Following the steps of Sections~\ref{sec:cancel-diag}
to~\ref{sec:M0x-diag} leads to a cubic
equation for $S(x)$ 
 which involves two
additional unknown series, namely $S_1$ and 
$$
F_0:=[x^0]\big(\Delta(x)S(x)S(\bx)\big)-3S(-1),
$$
with $\Delta(x)=1-4t(1+x)(1+\bx)$. This equation (which
is the counterpart of~\eqref{eqcatS-diag}) reads
$\Pol(S(x),S_1,F_0,t,x)=0$ with
\begin{multline*}
 \Pol(x_0,x_1,x_2,t,x)= \Delta(x) \left( x_0^3-\frac 3{x+1}x_0^2+ \frac{2-x}{x+1}x_0\right)\\
-\left( 16t^2x_1-x_2+t^2\bx +t^2x\right) x_0+7t^2x_1-x_2+t^2x+
\frac{x_2+2t^2x_1}{x+1}.
\end{multline*}
We apply to this equation the generalized quadratic method of
Section~\ref{sec:quad}. Two series, denoted $X_0$ and $X_1$, cancel
$\Pol_{x_0}$, and their first coefficients are:
$$
X_0=2-\frac{21}2 t^2-\frac{117}{8} t^4 + O(t^6), \qquad 
X_1= \frac 9 2 t^2 + \frac {261}8 t^4 + \frac{5067}{16}+ O(t^8).
$$
These coefficients do not suggest any obvious values for these series. To obtain
equations satisfied by $S_1$ and $F_0$, one has to actually work with
the system of 6 equations 
\begin{align*}
  \Pol(S(X_i), S_1,F_0,t,X_i)&=0,\\
\Pol_{x_0}(S(X_i), S_1,F_0,t,X_i)&=0,\\
 \Pol_x(S(X_i),S_1,F_0,t,X_i)&=0,
\end{align*}
 as explained in~\cite[Sec.~9]{mbm-jehanne}. The
most effective way seems to use Theorem~14 from~\cite{mbm-jehanne},
which says that the discriminant of $\Pol(x_0,  S_1,F_0,t,x)$ with respect
to $x_0$ admits $X_0$ and $X_1$ as \emm double, roots. Up to a
denominator and a factor 
$\Delta(x)$ (which does not vanish at $X_0$ nor $X_1$), this
discriminant is a polynomial of degree 4 in
$s:=x+\bx$, with coefficients in $\qs(t,S_1,F_0)$. This polynomial in $s$ has two double roots, namely $X_0+1/X_0$
and $X_1+1/X_1$, and thus it must be the square of a polynomial in $s$
of degree 2. This gives us two conditions on $S_1$ and $F_0$, from
which we obtain equations of degree 4 for each of these two series. As
before, they can be expressed rationally in terms of the series $T$:
$$
S_1= \frac{(T-1)(21-6T+T^2)}{(T+3)^3},
$$
$$
F_0=\frac{(T-1)(5T^3-11T^2+135T-33)}{128\, T^3}.
$$
From there one finds that $S(x)=t\sqrt x N(0, \sqrt x)$ (a series of
$\qs[x][[t^2]]$) has degree 12 
over $\qs(x,t)$, and is cubic over $\qs(T,x)$. It can be expressed
rationally in terms of $T$ and $V$.

For the series $R(x)=t^2/\sqrt x\, N(\sqrt x, 0)$, the degree is 24 over $\qs(x,t)$, and 3
over $\qs(Z,x)$.  Introducing $V$ factors the equation and gives
a rational expression of $R(x)$ in terms of $Z$ and $V$.

\section{Square lattice walks in a $\boldsymbol{135\degree }$ wedge}
\label{sec:gessel}
We now return to Ira Gessel's ex-conjecture~\eqref{conj:gessel} about square lattice walks starting and ending at $(0,0)$ and remaining in the
(convex) cone $\{(i,j): i+j\ge 0 \hbox{ and } j\ge 0\}$
(Figure~\ref{fig:gessel}, left).
More generally, let us denote by $g_{i,j}(n)$ the number of $n$-step walks in
this cone, starting at $(0,0)$ and ending at $(i,j)$. A step by step
construction of Gessel's walks gives
$$
\left(1-t(x+\bx+y+\by) \right) G(x,y)=1-t\by G(x,0) -t(\bx+\by)
G^\Delta(\bx y)+t\by G_{0,0},
$$
where 
$$
G^\Delta(x):= \sum_{j,n\ge 0} g_{-j,j}(n) x^j t^n
$$
counts walks ending on the diagonal $i+j=0$. Hence it suffices to
determine the series $G(x,0)$ and $G^\Delta(x)$.

\begin{figure}[b!]
\begin{center}
{\scalebox{1}{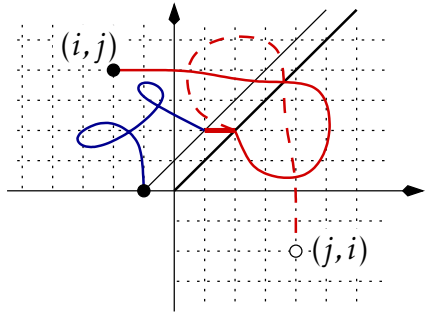}}
\end{center}
\caption{The reflection principle: a walk from $(-1,0)$ to $(i,j)$ that crosses the line $y=x+1$
  can be transformed bijectively into a walk ending at $(j,i)$.} 
\label{fig:refl}
\end{figure}

\medskip
One of our motivations for studying walks in a three-quadrant cone was
to attack the enumeration of Gessel's walks by the reflection
principle. Indeed, let us denote by $c_{i,j}(n)$
the number of $n$-step  walks going from $(-1,0)$ to $(i,j)$ on the
square lattice, and avoiding the negative quadrant. The \gf\ of these
numbers is given in Theorem~\ref{thm:square-asym}. Then the reflection principle
(Figure~\ref{fig:refl}) implies that, for $j \ge 0$ and $i<j$,
$$
c_{i,j}(n)-c_{j,i}(n)= g_{-i-1,j}(n).
$$

In particular, the case $j=0$ allows us to compute the specialization
$G(x,0)$ in terms of the series $L$ and $B$ of Theorem~\ref{thm:square-asym}:
$$
G(x,0)= L(x,0)-B(0,x).
$$
Using the rational expressions~\eqref{Lxt0} and~\eqref{B0xt} of $L(xt,0)$ and $B(0,xt)$
in terms of  $Z$ and $U$, we  recover the parametrized
expression of $G(xt,0)$ given in~\cite{BoKa09,mbm-gessel}.

\begin{figure}[hbt]
\begin{center}
{\scalebox{1}{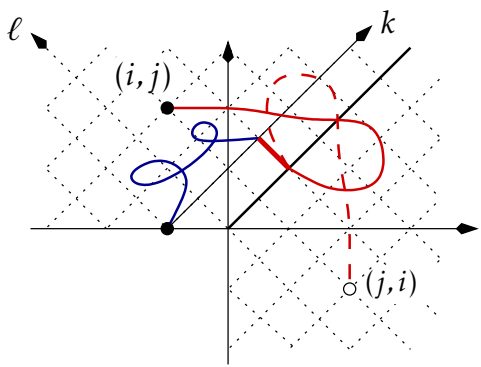}}
\end{center}
\caption{Second application of the reflection principle.} 
\label{fig:refl-diag}
\end{figure}

In order to determine the series $G^\Delta$, we have to extract from
$C(x,y)$  the
quasi-diagonal terms $c_{i,i+1}(n)$ and $c_{i+1,i}(n)$. To avoid
 this extraction, we can use instead our results on
the \emm diagonal, square lattice obtained in
Section~\ref{sec:diag-asym}. Indeed, let us denote by $\tilde c_{i,j}(n)$
the number of $n$-step walks going from $(-2,0)$ to $(i,j)$ on the
diagonal square lattice and avoiding the negative quadrant. The \gf\
of these numbers is given by Theorem~\ref{thm:diag-asym} (we use 
the tilde because we are mixing results for the square
lattice and the diagonal square lattice). The reflection principle
(Figure~\ref{fig:refl-diag}) now gives, for $j\ge 0$ and $i<j$:
$$
\tilde c_{i,j}(n)-\tilde c_{j,i}(n)= g_{k, \ell}(n)
$$
with $k=\frac{i+j}2+1$ and $\ell= \frac{j-i}2-1$.
In particular, the case $j=0$ gives us the value
$$
G^\Delta(x)= \frac 1{\sqrt{x}}\left(\tilde L(\sqrt x,0) -\ \tilde
B(0,\sqrt x)\right),
$$
where $\tilde L$ and $\tilde B$ are the series denoted $L$ and $B$ in
Theorem~\ref{thm:diag-asym}. The series $\sqrt{\bx}\tilde L(\sqrt
x,0)$ and $\sqrt{\bx} \tilde
B(0,\sqrt x)$ have rational expressions
in terms of $Z$ and $V$ (see~\eqref{Lx0-diag} and~\eqref{B0x-diag}), and we thus
recover  the parametrized
expression of $G^\Delta(x)$ given in~\cite{BoKa09,mbm-gessel}.

This solution of Gessel's model is only short because we have spent
much effort solving three-quadrant problems. The self-contained proofs
of~\cite{mbm-gessel,BeBoRa16} remain more direct.

\section{Questions, perspectives}


\subsection{About the present paper}
The first obvious problem raised by this paper is finding more
combinatorial proofs of our results. Since these results
include a solution to Gessel's famously difficult problem, this is
not likely to be easy. However, at least one question that arises from the
first  step of our approach should be easier.

Consider square lattice walks starting from $(0,0)$ and confined to
$\mathcal C$ (Theorem~\ref{thm:square}). The first two equations in
this theorem, namely~\eqref{C-split} and~\eqref{PM}, come at once by forming the orbit equation of $C(x,y)$,
and they imply that, for $i, j \ge 0$,
$$
C_{i,j}=Q_{i,j}+ C_{-i-2,j}+ C_{i,-j-2}.
$$
Given that forming the orbit equation is essentially taking reflections in
the coordinate axes, is there a simple explanation for this identity?
Note that it holds verbatim for walks starting at $(0,0)$ on the
diagonal square lattice. For square 
lattice walks starting at $(-1,0)$, the term in $Q$ disappears, leaving
$$
C_{i,j}= C_{-i-2,j}+ C_{i,-j-2}.
$$
For walks on the diagonal square lattice starting from $(-2,0)$, the
term in $Q$ remains, but its sign changes:
$$
C_{i,j}=-Q_{i,j}+ C_{-i-2,j}+ C_{i,-j-2}.
$$

\medskip
Another result that may be studied \emm per se, is the fact that
square lattice walks confined to $\mathcal C$, starting at $(0,0)$ and
ending at $(-1,0)$ are equinumerous with walks on the \emm diagonal,
square lattice, confined to $\mathcal C$ and joining $(0,0)$ to
$(-1,1)$ (see~\eqref{coincidence1} and~\eqref{coincidence2}).

\subsection{Perspectives}
As mentioned in the introduction, we hope that this paper will be the
starting point of a systematic study of walks with small steps
confined to $\mathcal C$, analogous to what has been achieved in the
past decade for walks confined to the first quadrant~$\mathcal Q$. By
small steps, we mean steps taken from $\{-1,0,1\}^2$. 

\begin{figure}[htb]
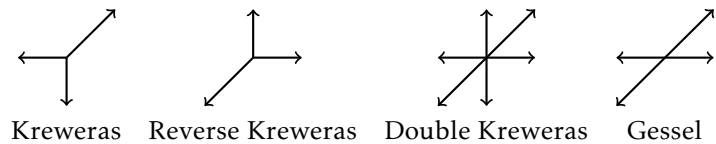

  \centering
  \begin{tabular}{cccc}
     $ \diagr{NE,S,W}$  & $  \diagr{SW,N,E}$
  & $\diagr{N,S,E,W,NE,SW} $ &
  $ \diagr{E,W,SW,NE}$
\\
Kreweras & Reverse Kreweras& Double Kreweras& Gessel
  \end{tabular}
   \caption{The four algebraic quadrant models.}
  \label{fig:alg}
\end{figure}

Let us recall some of the quadrant results: given a set $\cS$
of small steps, the \gf\ $Q(x,y)$ that counts walks starting from
$(0,0)$, confined to $\mathcal Q$ and taking their steps in $\cS$ is
D-finite if and only if a certain group of rational transformations is
finite. This happens for 23 inherently different step sets, among
which exactly 4 even lead to an algebraic \gf\ (Figure~\ref{fig:alg}). There remain 56
inherently different non-D-finite step sets, among which 5 are called
\emm singular,: this means that all their elements $(i,j)$ satisfy
$i+j\ge 0$. Two generic approaches prove the non-D-finiteness of the
51 non-singular models~\cite{KuRa12,BoRaSa12}, and the remaining 5 are proved non-D-finite in
a more \emm ad hoc, way~\cite{MiRe09,MeMi13}. 

Now what happens for walks with small steps confined to the three-quadrant cone $\cC$?
\begin{itemize}
\item One can check that all models that are trivial or simple when
  counting walks   confined to the first quadrant (and have a rational
  or algebraic \gf\ for elementary reasons~\cite{BoMi10}) are still trivial or
  simple when counting walks avoiding the negative quadrant. This
  leaves us, as in the quadrant problem, with 79 inherently different models.
\item The case of singular step sets is particularly simple: all walks
  formed of such steps remain in the half-plane $i+j\ge 0$, and \emm a
  fortiori, in $\cC$. Hence the associated \gf\ is rational, equal to
  $\left(1-t\sum_{(i,j)\in \cS}x^iy^j\right)^{-1}$. A simple start!
\item Could it be that for any step set associated with a finite
  group, the \gf\ $C(x,y)$ is D-finite? and, maybe, differs from a
  simple D-finite series related to $Q(x,y)$ by an algebraic series?
\item In particular, could it be that for the four step sets of
  Figure~\ref{fig:alg},  for which $Q(x,y)$ is known to be algebraic, $C(x,y)$
  is also algebraic? We could not resist trying a bit of guessing on
  these models, and algebraicity seems very plausible. At least, we
  have guessed in each case an algebraic equation for the series
  $C_{0,0}$ that counts walks starting and ending at $(0,0)$. 
 The  degree is, from left to right, 6, 6, 16, 24, which should be
 compared to the values 3, 3, 4, 8 obtained for quadrant walks. 
\item To what extent can the approach of this paper be adapted to
  other step sets associated with a finite group? A first candidate
  would be the set of all eight small steps, which has the same
  symmetries as the models studied here, and is quite likely to be solvable by
  the same approach.
\item Could it be that for non-singular step sets associated with an
  infinite group, the series $C(x,y)$ is non-D-finite? Can this be
  proved using asymptotic enumeration, as
  has been done for quadrant walks~\cite{BoRaSa12} using the results of~\cite{denisov-wachtel}?
  (This question has been answered positively by Mustapha after this
  paper appeared on  arXiv~\cite{mustapha}.)
\item  Can the powerful analytic approach of~\cite{raschel-unified} be
  adapted to walks avoiding a quadrant? This approach was the first to
  yield non-D-finiteness results for the 51 models non-singular with
  an infinite group~\cite{KuRa12}. 
\end{itemize}

To finish, let us mention that the reflection principle relates
several three-quadrant models to quadrant models, as exemplified in
Section~\ref{sec:gessel}. More precisely, counting walks with Kreweras
steps in $\cC$ gives a solution of walks with \emm reverse, Kreweras
steps in $\mathcal Q$, and vice-versa. Similarly, counting walks with double Kreweras
steps in $\cC$ also solves walks with double Kreweras steps in
$\mathcal Q$.

\bigskip\bigskip
\noindent{\bf Acknowledgements.} Back in 2007, the author  had many interesting
discussions with Ira Gessel about the closed forms~\eqref{expressions-diag}, which were at that time conjectures. Ira produced more
conjectures of this type for other endpoints. 

 \bibliographystyle{plain}
\bibliography{qdp}


\appendix
\section{Parametrized  expressions}
Our parametrizing series $T$, $Z$, $U$ and $V$ are defined as follows.
First, $T$ is the unique series in $t$ with
constant term $1$ satisfying:
$$
T= 1+256\, t^2 \frac{T^3}{(T+3)^3},
$$
and  $Z=\sqrt T$. We have:
$$
T=1+4\,{t}^{2}+36\,{t}^{4}+396\,{t}^{6}+4788\,{t}^{8}+O \left( {t}^{10}
 \right) ,
$$
$$
Z=1+2\,{t}^{2}+16\,{t}^{4}+166\,{t}^{6}+1934\,{t}^{8}+O \left( {t}^{10}
 \right) 
.
$$
In fact, $Z$ is the sum of two hypergeometric series:
$$
Z=\sum_{n\ge 0} {16^n}\left(2\, \frac{(-1/2)_{n}(1/6)_{n}}{(1)_n(1/3)_n}
-
\frac{(-1/2)_{n}(5/6)_{n}}{(1)_n(2/3)_n}
\right) t^{2n}.
$$

Then,  $U$ is the only power series in $t$ with constant term $1$ (and
coefficients in $\qs[x]$) satisfying 
$$
16T^2(U^2-T)=x(U+UT-2T)(U^2-9T+8TU+T^2-TU^2).
$$

Finally,  $V$ is the only series in $t$ with constant term $0$ (and
coefficients in $\qs[y]$)
satisfying 
$$
1-T+3V+VT= yV^2(3+V+T-VT).
$$ 
We have
$$
U=1+2\,{t}^{2}+16\,{t}^{4}+ \left( 166+2\,x \right) {t}^{6}+ \left( 2\,
{x}^{2}+40\,x+1934 \right) {t}^{8}+O \left( {t}^{10} \right) ,
$$
$$
V={t}^{2}+ \left( 8+x \right) {t}^{4}+ \left( 2\,{x}^{2}+16\,x+82
 \right) {t}^{6}+ \left( 5\,{x}^{3}+48\,{x}^{2}+227\,x+944 \right) {t}
^{8}+O \left( {t}^{10} \right) 
.
$$

\subsection{Walks starting at $\boldsymbol{(0,0)}$ on the square
  lattice} 
\label{sec:param-sq}
The \gf\ of walks ending on the negative $x$-axis (resp. at
abscissa~$-1$) is $\bx M(\bx,0)- \bx^2Q(\bx,0)/3$ (resp. $\bx M(0, y)$) where 
$Q(x,y)$ is given by~\eqref{Q-expr} and
\beq\label{Mxt0}
tM(xt,0)=\frac{N_1(Z,U)}{3T(Z-1)(T+3)^3(U+Z)^2(U^2-9T+8TU+T^2-TU^2)}
\eeq
and
\beq\label{M0xt}
tM(0,xt)=\frac{(TU-2T+U)^2N_2(T,U)}{T(T+3)^3D_2(T,U)},
\eeq
with
\begin{multline*}
N_1(z,u)=- \left( {z}^{2}+1 \right) ^{2} \left( z+1 \right) ^{3} \left( z-1
 \right) ^{4}{u}^{8}\\
-2\,z \left( {z}^{2}+1 \right)  \left( {z}^{4}-10
\,{z}^{3}-14\,z-1 \right)  \left( z+1 \right) ^{2} \left( z-1 \right) 
^{3}{u}^{7}\\
+4\,{z}^{2} \left( z+1 \right)  \left( 10\,{z}^{7}-35\,{z}^
{6}+14\,{z}^{5}-115\,{z}^{4}-10\,{z}^{3}-57\,{z}^{2}-14\,z+15 \right) 
 \left( z-1 \right) ^{2}{u}^{6}\\
+2\,{z}^{3} \left( z-1 \right)  \left( 
{z}^{10}+14\,{z}^{9}-77\,{z}^{8}+252\,{z}^{7}-66\,{z}^{6}+1224\,{z}^{5
}+106\,{z}^{4}+68\,{z}^{3}\right. \\ \left. +33\,{z}^{2}-534\,z+3 \right) {u}^{5}+2\,{z}
^{4} \left( z-1 \right)  \left( {z}^{10}+8\,{z}^{9}-115\,{z}^{8}-400\,
{z}^{7}-1154\,{z}^{6}\right. \\ \left.-1728\,{z}^{5}-5890\,{z}^{4}-1520\,{z}^{3}-2607\,
{z}^{2}-456\,z+549 \right) {u}^{4}+2\,{z}^{5} \left( {z}^{11}-11\,{z}^
{10}\right. \\ \left.-207\,{z}^{9}+149\,{z}^{8}+2946\,{z}^{7}+2202\,{z}^{6}+8506\,{z}^{
5}+9266\,{z}^{4}-5571\,{z}^{3}+4017\,{z}^{2}\right. \\ \left.-3627\,z-1287 \right) {u}^
{3}-4\,{z}^{6} \left( 14\,{z}^{10}-{z}^{9}-465\,{z}^{8}-684\,{z}^{7}+
3704\,{z}^{6}+2034\,{z}^{5}\right. \\ \left.+11274\,{z}^{4}+6756\,{z}^{3}-702\,{z}^{2}+
3159\,z-513 \right) {u}^{2}-2\,{z}^{7} \left( {z}^{11}+13\,{z}^{10}-
115\,{z}^{9}\right. \\ \left.-95\,{z}^{8}+1346\,{z}^{7}+3722\,{z}^{6}-8334\,{z}^{5}-
4470\,{z}^{4}-24291\,{z}^{3}-12663\,{z}^{2}\right. \\ \left.-351\,z-3915 \right) u-{z}^
{8} \left( {z}^{11}-{z}^{10}-37\,{z}^{9}+293\,{z}^{8}+382\,{z}^{7}-894
\,{z}^{6}-4614\,{z}^{5}\right. \\ \left.+7686\,{z}^{4}+5409\,{z}^{3}+17631\,{z}^{2}+
9099\,z-2187 \right) ,
\end{multline*}
\begin{multline*}
 N_2(t,u)={t}^{5}+ \left( -2\,{u}^{2}+32\,u-73 \right) {
t}^{4}+ \left( {u}^{4}-16\,{u}^{3}+90\,{u}^
{2}-96\,u-177 \right) {t}^{3}\\+ \left( -{u}^{4}+82
\,{u}^{2}-192\,u-135 \right) {t}^{2}-{u}^{
2} \left( {u}^{2}-16\,u+42 \right) t+{u}^{
4},
\end{multline*}
and
$$
  D_2(t,u)=-{t}^{4}+2\, \left( u-2 \right)  \left( u-6
 \right) {t}^{3}- \left( u-3 \right)  \left( {u}^
{3}+3\,{u}^{2}-15\,u+3 \right) {t}^{2} +6\,t\,{u}^{2}+{u}^{4}.
$$
%

\subsection{Walks starting at $\boldsymbol{(0,0)}$ on the diagonal
  square lattice}
\label{sec:param-diag}
The \gf\ for walks ending on the negative $x$-axis (resp. at
abscissa~$-1$) is $\bx M(\bx,0)- \bx^2Q(\bx,0)/3$ (resp. $\bx M(0,
y)$) where $Q(x,y)$ is given by~\eqref{Q-diag} and
\beq\label{Mx0-diag}
\frac {t^2} {\sqrt x}\, M( \sqrt x,0)
=\frac{V(VT-T-V-3)N(Z,V)}{48Z^3(V+1)^2(VZ-V-Z-1)^2},
\eeq
\beq\label{M0x-diag}
t\sqrt x\, M(0, \sqrt x)
=
\frac{(TV-T+3V+1)(-TV^2-2TV+V^2+T+2V+3)}
{2(V+1)(T^2V^2-2T^2V+2TV^2+T^2-3V^2+2T-6 V-3)}
\eeq
with
\begin{multline*}
  N(z,v)=- \left( {z}^{2}+3 \right)  \left( z-1 \right) ^{2}{v}^{3}+ \left( z-1
 \right)  \left( 3\,{z}^{3}+9\,{z}^{2}+z+11 \right) {v}^{2}\\+ \left( z+
1 \right)  \left( 3\,{z}^{3}-9\,{z}^{2}+z-11 \right) v- \left( {z}^{2}
+3 \right)  \left( z+1 \right) ^{2}.
\end{multline*}

\subsection{Walks starting at $\boldsymbol{(-1,0)}$ on the square
  lattice}
\label{sec:param-sq-asym}
The \gf\ for walks ending on the negative $x$-axis (resp. at
abscissa~$-1$, on the negative $y$-axis, at ordinate $-1$) is $\bx
L(\bx,0)$ 
(resp. $\bx L(0, y)$, $\by B(0,\by)$, $\by B(x,0)$) where 
\beq\label{Lxt0}
L(xt,0)=\frac{256 Z^4(UZ^2-2Z^2+U)}{(UZ+Z^2+U-3Z)(U+Z)(Z-1)(Z^2+3)^3},
\eeq
$$
xL(0,xt)=\frac{256Z^6(U-Z)(UZ-Z^2-U-3Z)(U^2Z^2-Z^4-4UZ^2+U^2+3Z^2)}
{D(Z,U)(Z-1)(Z^2+3)^3},
$$
$$
xB(xt,0)=\frac{512Z^6(UZ^2-2Z^2+U)(U-Z)N_1(Z,U)}{(Z^2-1)(Z^2+3)^3D(Z,U)},
$$
and
\beq\label{B0xt}
B(0,xt)=\frac{16Z^2(UZ^2-2Z^2+U) N_2(Z,U)}{(U+Z)(1-Z^2)(Z^2+3)^3(UZ-Z^2-U-3Z)}
\eeq
with
$$
D(z,u)={z}^{8}-2\, \left( u-2 \right)  \left( u-6 \right) {z}^{6}+ \left( u-3
 \right)  \left( u^{3}+3\,u^{2}-15\,u+3 \right) {z}^{4}-6\,u^{2}
{z}^{2}-u^{4}
,
$$
$$
N_1={z}^{4}+ 2\left( u-3 \right) {z}^{3}+ \left( u+1 \right)  \left( u-3
 \right) {z}^{2}-u^{2},
$$
$$
N_2(z,u)=\left( z+1 \right)  \left( z-1 \right) ^{2}u ^{2}-4\,z \left( 3\,z+1
 \right)  \left( z-1 \right) u -{z}^{2} \left( {z}^{3}+3\,{z}^{2}-25\,z
-11 \right)
.
$$
\subsection{Walks starting at $\boldsymbol{(-2,0)}$ on the diagonal
  square lattice}
\label{sec:param-diag-asym}
The \gf\ for walks ending on the negative $x$-axis (resp. at
abscissa~$-1$, on the negative $y$-axis, at ordinate $-1$) is $\bx
L(\bx,0) + \bx^2Q(\bx,0)/3$
(resp. $\bx L(0, y)$, $\by B(0,\by)+ \by^2Q(\by,0)/3$, $\by B(x,0)$)  where $Q(x,y)$ is given by~\eqref{Q-diag} and
\beq\label{Lx0-diag}
\frac 1 {\sqrt x} L( \sqrt x,0)
=\frac{32VZ^3(VT-T-V-3)N_1(Z,V)}{3(1+V)^2(T-1)(T+3)^3(VZ-V-Z-1)^2},
\eeq
$$
t\sqrt x L(0, \sqrt x)= \frac{(TV-T+3V+1)(V-1)(TV-T-V-3)}
{2(V+1)(T^2V^2-2T^2V+2TV^2+T^2-3V^2+2T-6 V-3)},
$$
$$
t {\sqrt x} B( \sqrt x,0)
=\frac{V(TV-T+3V+1)(2-TV+V)}
{(V+1)(T^2V^2-2T^2V+2TV^2+T^2-3V^2+2T-6 V-3)},
$$
\beq\label{B0x-diag}
\frac 1 {\sqrt x} B(0, \sqrt x)
=\frac{16VZ^3(VT-T-V-3)N_2(Z,V)}{3(V+1)^2(T-1)(T+3)^3(VZ-V-Z-1)^2}
\eeq
with
\begin{multline*}
  N_1(z,v)=- \left( {z}^{2}+3 \right)  \left( z-1 \right) ^{2}{v}^{3}+4 \left( 
z+2 \right)  \left( z-1 \right) {v}^{2}
+4 \left( z+1 \right) 
 \left( z-2 \right) v- \left( {z}^{2}+3 \right)  \left( z+1 \right) ^{
2}
\end{multline*}
and
\begin{multline*}
  N_2(z,v)=\left( z^{2}+3 \right)  \left( z-1 \right) ^{2}v^{3}+ \left( 3{
z}^{2}+12z+5 \right)  \left( z-1 \right) ^{2}v ^{2}\\
+ \left( 3z^
{2}-12  z+5 \right)  \left( z+1 \right) ^{2}v + \left( z^{2}+3
 \right)  \left( z+1 \right) ^{2}.
\end{multline*}

\end{document}